# A generalized inverse for graphs with absorption

Karly Jacobsen [*]   Joseph Tien [*†]

November 4, 2016

**Abstract**

We consider weighted, directed graphs with a notion of absorption on the vertices, related to absorbing random walks on graphs. We define a generalized inverse of the graph Laplacian, called the absorption inverse, that reflects both the graph structure as well as the absorption rates on the vertices. Properties of this generalized inverse are presented, including a matrix forest theorem relating this generalized inverse to spanning forests of a related graph, as well as relationships between the absorption inverse and the fundamental matrix of the absorbing random walk. Applications of the absorption inverse for describing the structure of graphs with absorption are presented, including a directed distance metric, spectral partitioning algorithm, and centrality measure.

## 1 Introduction

This paper concerns a generalized inverse of the graph Laplacian that is related to transient random walks on graphs. Let $G$ denote a strongly connected, weighted, directed graph on the vertex set $V = \{1, 2, \ldots, n\}$. Let $A \in \mathbb{R}^{n \times n}$ be the weighted adjacency matrix with $ij$th entry $a_{ij}$ corresponding to the weight of the edge from $j$ to $i$. The outdegree (i.e. sum of outgoing weights) from vertex $i$ is given by $w_i = \sum_{j=1}^n a_{ji}$. The graph Laplacian matrix is given by $L = W - A$ where $W = \text{diag}\{w\}$, the diagonal matrix with $W_{ii} = w_i$. Both the group inverse, $L^{\#}$, and the Moore-Penrose pseudoinverse, $L^{\dagger}$, of the graph Laplacian have been shown to encode valuable information about the graph as well as the regular Markov process generated by $L$ including hitting times, measures of centrality, distance functions on trees, and more [1–5].

Here we study a generalization of the group inverse that similarly captures valuable information about transient random walks on the graph, where each vertex is a transient state with transition rate $d_i > 0$ to an absorbing state. We will refer to the $d_i$ as *absorption rates*. Let $d = (d_1, \ldots, d_n)^T$ denote the (column) vector of absorption rates, and let $D = \text{diag}\{d\}$. Formally, we define a *graph with absorption* as the pair $(G, d)$.

The object of interest for this paper is the matrix $X$ which satisfies the two conditions listed below:

$$XLy = y \text{ for } y \in N_{1,0}, \tag{1}$$
$$Xy = 0 \text{ for } y \in R_{1,0}, \tag{2}$$

where

$$N_{1,0} = \{x \in \mathbb{R}^n : Dx \in \text{range } L\}, \tag{3}$$
$$R_{1,0} = \{Dx : x \in \ker L\}. \tag{4}$$

That is, $X$ sends $R_{1,0}$ to 0 and acts as the left inverse of $L$ on $N_{1,0}$.

**Definition 1.** *Let $L$ be the graph Laplacian of a strongly connected graph with absorption $(G, d)$, with $d > 0$. A matrix $X$ satisfying (1)–(2) is said to be an absorption inverse of $L$ with respect to $d$.*

---
[*]Mathematical Biosciences Institute, The Ohio State University, Columbus, OH 43210, USA
[†]Department of Mathematics, The Ohio State University, Columbus,OH 43210, USA

We will denote the absorption inverse of $L$ with respect to $d$ by $L^d$. As the name implies, $L^d$ is indeed a generalized inverse of $L$. Specifically, Tien et al. [6] show that $L^d$ is a $\{1,2\}$-inverse of $L$. Derivation of $L^d$ in [6] is through a Laurent series expansion, using the elegant results of Langenhop [7]. Indeed, Langenhop shows that conditions (1) and (2) arise naturally in the context of Laurent series expansion for the inverse of a nearly singular matrix. This motivates our definition of $L^d$.

Our interest in studying $L^d$ is due to its connections with transient random walks on graphs, and in particular with how the absorption rates $d_i$ influence how we think of distance and other measures of structure on a graph with absorption. To give a simple illustration of this, consider an undirected cycle with $n$ vertices, labeled in ascending order clockwise. Topologically the graph is a circle, and a random walker on the graph starting at 1 can reach $n-1$ by passing through $n$. Now suppose that we introduce absorption to the graph with a very large absorption rate at $n$. This effectively puts a barrier between 1 and $n-1$, as a walker starting at 1 that arrives at $n$ will have a very high probability of absorption. It will thus be easier to reach $n-1$ by first traversing all other vertices on the graph clockwise between 1 and $n-1$. The distance between 1 and $n-1$ has effectively increased due to the absorption at $n$, and indeed the graph has essentially changed from a circle to a path. These changes in structure are not captured by traditional metrics based solely upon the arc weights of the graph. Instead, metrics that take into account both graph structure and absorption are needed. Tien et al. [6] suggest that $L^d$ may provide natural ways to do so. Graphs with absorption arise naturally in many applications including birth-death processes, epidemiology, and ecology, and thus $L^d$ is of applied interest. For example, disease spread on a community network depends both on the network structure and the rate of pathogen decay (i.e. absorption) in each community. The absorption inverse also appears in the series expansion for the fundamental matrix of nearly decomposable Markov chains [8–11], and an expression relating $L^d$ to the Moore-Penrose inverse is given in [10].

The purpose of this paper is to establish some basic properties of $L^d$, and illustrate how $L^d$ can be used to examine the structure of graphs with absorption. Background material on graph theory, generalized inverses, and previous results on $L^d$ are given in Section 2 for the convenience of the reader. Our results include simple formulas relating $L^d$ to the group and Moore-Penrose inverses (Section 3.1) as well as formulas relating $L^d$ to $L$ that hold outside the radius of convergence of the Laurent series (Section 3.2). In particular, we give a characterization of $L^d$, that holds for all $d > 0$, in terms of the spanning forests of a related graph (Section 3.3). We also introduce a new distance metric based on $L^d$ that takes differential absorption rates into account, together with an extension of PageRank [12] giving vertex centrality for graphs with absorption (Sections 4.1−4.2). Connections with spectral partitioning methods based on $L^d$ are also discussed (Section 4.3). Simple examples illustrating these structural measures and algorithms are given in Section 4.4. Concluding remarks are given in Section 5, and explicit formulas for $L^d$ for several common graph motifs are given in A.

## 2 Preliminaries

We present here some definitions and basic facts from graph theory and generalized inverses, together with previous results about $L^d$.

### 2.1 Graph theory

A directed graph is *strongly connected* if for any ordered pair of distinct vertices, there exists a path from one to the other. A directed graph is strongly connected if and only if the weighted adjacency matrix $A$ is irreducible. A graph is said to be *balanced* (also called vertex balanced [13], or Eulerian [14]) if the sum of incoming edge weights equals the sum of outgoing edge weights for each vertex [6]. Our work in Section 4 will be restricted to the case of balanced graphs.

Many properties of a graph can be examined through study of an associated graph Laplacian. Several different versions of the Laplacian matrix exist in the literature. The version we use here $L = W - A$ is the "unnormalized Laplacian", in the terminology of [2]. The dimension of $\ker L$ is the number of connected components of $G$, and hence is one for a strongly connected graph. The matrix tree theorem implies that there exists $u = (u_1, \ldots, u_n)^T$, a positive basis for $\ker L$ with $\sum_{i=1}^{n} u_i = 1$ and $u_i$ proportional to the sum of the weights of all spanning trees rooted at vertex $i$ ([6], Appendix A). For balanced $G$, $u_i = 1/n$ for all $i$. Note also that $\mathbf{1}^T L = 0$ where $\mathbf{1}$ denotes the $n-$dimensional column vector of ones.

Table 1: Summary of notation

| Notation | Description |
|---|---|
| $G$ | graph on vertex set $\{1, 2, ..., n\}$ |
| $d$ | positive vector of absorption rates |
| $\bar{d}$ | $d^T u$ |
| $D$ | diag$\{d\}$, diagonal matrix of absorption rates |
| $n$ | number of nodes in $G$ |
| $\tilde{G}$ | absorption-scaled graph |
| $A$ ($\tilde{A}$) | adjacency matrix of graph $G$ ($\tilde{G}$) |
| $a_{ij}$ ($\tilde{a}_{ij}$) | weight of edge in $G$ ($\tilde{G}$) from $j$ to $i$ |
| $w_i$ ($\tilde{w}_i$) | outdegree of vertex $i$ in $G$ ($\tilde{G}$) |
| $W$ ($\tilde{W}$) | diagonal matrix of outdegrees for graph $G$ ($\tilde{G}$) |
| $L$ ($\tilde{L}$) | graph Laplacian of $G$ ($\tilde{G}$) |
| $L^\dagger$ | Moore-Penrose pseudoinverse of $L$ |
| $L^\#$ | group inverse of $L$ |
| $L^d$ | absorption inverse of $L$ |
| $N_{1,0}$ | $\{x \in \mathbb{R}^n : Dx \in \text{range } L\}$ |
| $R_{1,0}$ | $\{Dx : x \in \ker L\}$ |
| $u$ ($\tilde{u}$) | positive basis for $\ker L$ ($\ker \tilde{L}$) |
| $\omega(\Gamma)$ | weight of subgraph $\Gamma$ or nonempty set of subgraphs |
| $\mathcal{F}_k$ | set of all in-forests of $G$ with $k$ arcs |
| $\mathcal{F}_k^{i \to j}$ | in-forests in $\mathcal{F}_k$ in which $i$ belongs to a tree converging to $j$ |
| $Q_k$ | matrix of in-forests with $k$ arcs |
| $\sigma_k$ | weight of all in-forests with $k$ arcs, i.e. $\omega(\mathcal{F}_k)$ |
| $J_k$ | $Q_k/\sigma_k$, normalized matrix of in-forests with $k$ arcs |
| $U$ | $u\mathbf{1}^T/\bar{d}$ |
| $P$ | probability transition matrix of jump process |
| $\pi$ | stationary distribution of jump process |
| $\Pi$ | diag$\{\pi\}$ |
| $Z$ | fundamental matrix for regular Markov process |
| $\hat{M}$ | the bottleneck matrix of $L$ based at vertex $n$ |
| $R$ | absorption-scaled forest distance |

The Laplacian is also connected to spanning forests of $G$, for example through the all-minors tree theorem [15]. Chebotarev and Agaev study spanning trees and forests of a weighted, directed graph and their connections with the Laplacian matrix [16]. We present here some important definitions, adapting the notation of [16]. For a subgraph $\Gamma$ of $G$, let $\omega(\Gamma)$ denote the weight of $\Gamma$, given by the product of the weights of all its arcs. The weight of a nonempty set of subgraphs, $\mathcal{G}$, is given by $\omega(\mathcal{G}) = \sum_{\mathcal{H} \in \mathcal{G}} \omega(\mathcal{H})$. A *converging tree* is a directed graph that is weakly connected (i.e. its undirected version is connected) with one vertex, the "root", having outdegree 0 and all other vertices having a single outgoing edge. An *in-forest* of $G$ is a spanning directed graph in which all components are converging trees. Let $\mathcal{F}_k$ denote the set of all in-forests of $G$ with $k$ arcs (i.e. those with $n-k$ trees) and let $\mathcal{F}_k^{i \to j} \subset \mathcal{F}_k$ denote those in-forests in which vertex $i$ belongs to a tree converging to root $j$. Let $Q_k$ denote the *matrix of in-forests with $k$ arcs*, which has $ij$th entry given by $[Q_k]_{ij} = \omega(\mathcal{F}_k^{j \to i})$. Lastly, we denote by $\sigma_k$ the weight of all in-forests with $k$ arcs, i.e. $\sigma_k = \omega(\mathcal{F}_k)$.

A parametric matrix-forest theorem is presented in [16] that relates $(I + \tau L)^{-1}$, for any $\tau \in \mathbb{R}$, to the forest matrices of the graph:

**Theorem 1** (Theorem 3' in [16]). *For any $\tau \in \mathbb{R}$,*

$$(I + \tau L)^{-1} = \frac{1}{\sigma(\tau)} \left( Q_0 + \tau Q_1 + \ldots + \tau^{n-1} Q_{n-1} \right)$$

*where $\sigma(\tau) = \sum_{k=0}^{n-1} \sigma_k \tau^k$.*

Chebotarev and Agaev go on to express the forest matrices as polynomials in the Laplacian matrix and also to consider the relationship between forest matrices and generalized inverses of $L$. In particular, they show that the group inverse of $L$ can be expressed as

$$L^{\#} = (L + \tau \tilde{J})^{-1} - \tau^{-1} \tilde{J}, \qquad \tau \neq 0 \tag{5}$$

where $\tilde{J} = \sigma_{n-1}^{-1} Q_{n-1}$ is the eigenprojection at 0 of $L$, i.e. the projection on $\ker L$ along range $L$.

## 2.2 Generalized inverses

A matrix $X^+$ is called a $\{1,2\}$-inverse of a matrix $X$ if $X^+$ satisfies the following two conditions [17]:

$$XX^+X = X, \tag{6}$$
$$X^+XX^+ = X^+. \tag{7}$$

More generally, a $\{1\}$-inverse of $X$ satisfies (6). The *index* of a square matrix $X$ is defined as the smallest integer $k$ such that rank $X^k$ = rank $X^{k+1}$. The group inverse of a square matrix $X$ of index 1 is the unique matrix $X^{\#}$ that satisfies (6), (7) and the additional condition $XX^{\#} = X^{\#}X$. The Moore-Penrose pseudoinverse of $X$ is the unique matrix $X^{\dagger}$ that satisfies (6), (7) and the conditions $(XX^{\dagger})^T = XX^{\dagger}$ and $(X^{\dagger}X)^T = X^{\dagger}X$ where $^T$ denotes the conjugate transpose.

The group inverse $L^{\#}$ inverts $L$ on range $L$ and sends $\ker L$ to zero [17]. The generalized inverse $L^d$ can then be interpreted as a generalization of $L^{\#}$: $N_{1,0}$ and $R_{1,0}$ are distortions of range $L$ and $\ker L$, respectively, through the absorption rates $D$. As pointed out by [6], when the absorption rates $d_i$ are equal for all vertices, $N_{1,0} = $ range $L$, $R_{1,0} = \ker L$, and $L^d = L^{\#}$.

The connection between $L^d$ and the group inverse can also be seen by considering the resolvent of $L$. Ben-Israel [18] shows that as $z \to 0$,

$$\lim_{z \to 0} (L + zI)^{-1} L = L^{\#} L. \tag{8}$$

We likewise consider $(L + zD)^{-1}$ as a generalization of the resolvent of $L$. As shown by Langenhop [7] and Tien et al. [6], there exists a Laurent series with a simple pole given by:

$$(L + zD)^{-1} = \frac{1}{z} U + L^d + \sum_{k=1}^{\infty} (-zL^d D)^k L^d \tag{9}$$

where $\bar{d} = d^T u = \sum_{i=1}^{n} d_i u_i$ and $U \in \mathbb{R}^{n \times n}$ is given by $U = u\mathbf{1}^T / \bar{d}$. That is, the columns of $U$ are identical and proportional to $u$, the positive basis for $\ker L$. For $z > 0$, the convergence criterion for the series is $\rho(zL^d D) < 1$, where

$\rho$ denotes the spectral radius. A generalization of (8) follows directly from equation (9):

$$\lim_{z \to 0}(L + zD)^{-1}L = L^d L. \tag{10}$$

For $L \in \mathbb{R}^{n \times n}$ be of index one, $\mathbb{R}^n = \text{range } L \oplus \ker L$. Similarly, there are decompositions for $\mathbb{R}^n$ in terms of $N_{1,0}$ and $R_{1,0}$ [7]. These, and a few more results of Langenhop that will be needed, are stated below in Theorem 2.

**Theorem 2** (Langenhop [7]).

(a) $\mathbb{R}^n = N_{1,0} \oplus \ker L$.

(b) $\mathbb{R}^n = R_{1,0} \oplus \text{range } L$.

(c) $L$ is a bijection from $N_{1,0}$ onto $\text{range } L$.

(d) $L^d L = I - UD$ is a projection onto $N_{1,0}$, and $L^d$ sends $R_{1,0}$ to 0.

(e) $UD$ is a projection onto $\ker L$, which sends $N_{1,0}$ to 0.

We next prove a straightforward lemma that is analogous to Theorem 2d.

**Lemma 1.** *(a) $LL^d = I - DU$ is a projection onto $\text{range } L$, which sends $R_{1,0}$ to 0.*

*(b) $DU$ is a projection onto $R_{1,0}$, which sends $\text{range } L$ to 0.*

*Proof.* Let $x \in \text{range } L$ be given by $x = L\omega$. Then $LL^d x = LL^d L\omega = L\omega = x$. Likewise, $(I - DU)x = x - \frac{1}{d}Du\mathbf{1}^T L\omega = x$. Let $y \in R_{1,0}$. That is, $y = \beta Du$ for some $\beta$. Then $LL^d y = 0$ by Theorem 2d and $(I - DU)y = y - \frac{\beta}{d}Du\mathbf{1}^T Du = y - \beta Du = 0$. The claim *(a)*, and hence *(b)* as well, then follows from Theorem 2b. □

We conclude this section by providing the following formula for $L^d$, given by Tien et al. [6]:

**Proposition 1.** *(Tien et al [6]) Let $B : \mathbb{R}^n \to \mathbb{R}^n$ be the linear map defined by*

$$Bx = \begin{cases} x, & x \in \text{range } L \\ UDx, & x \in R_{1,0}. \end{cases}$$

*Then $L^d = (I - UD)(L - UD)^{-1} B$.*

## 2.3 Random walks on graphs

The Laplacian $L$ generates a regular random walk on $G$ with stationary distribution $u$, i.e. the infinitesimal generator for the random walk is given by $Q = -L$. The associated jump process (i.e. embedded Markov chain) corresponds to a random walk on the weighted, directed graph with transition matrix $P = AW^{-1}$. Correspondingly, the jump chain has a stationary distribution, $\pi$, such that $P\pi = \pi$. Let $\Pi = \text{diag}\{\pi\}$.

For a graph with absorption $(G, d)$, consider an associated Markov process on $n + 1$ states where the vertices $1, \ldots, n$ represent the transient states, $n+1$ an absorbing state, and transitions from the transient to absorbing state correspond to the absorption rates $d$. The infinitesimal generator, $Q_d \in \mathbb{R}^{(n+1) \times (n+1)}$, for this transient random walk is then

$$Q_d = \begin{bmatrix} -(L + D) & \mathbf{0} \\ d^T & 0 \end{bmatrix}.$$

where $\mathbf{0}$ is the $n$-dimensional column vector of zeros. The matrix $(L + D)^{-1}$ is called the *fundamental matrix of the absorbing process* [19]. The $ij$th entry, $(L+D)^{-1}_{ij}$, represents the expected time spent (i.e. *residence time*) in state $i$ before absorption occurs, given that the process started in state $j$ [20].

The group and Moore-Penrose inverses of $L$ are closely connected with random walks on $G$. For example, the entries of $L^\dagger$ are related to the *commute times* of the random walk [21], as discussed in Section 3.4. For a regular Markov process the $ij$th entry of $L^\#$ can be interpreted as the deviation from the expected time spent in vertex $i$, due to starting from vertex $j$ [5].

The absorption inverse $L^d$ is similarly related to fluctuations in the time spent in $i$ starting from $j$ for the *transient* random walk on $(G, d)$. We take $z = 1$ in (9) to get a series expansion for the fundamental matrix for the transient random walk:

$$(L + D)^{-1} = U + L^d + \sum_{k=1}^{\infty} (-L^d D)^k L^d, \tag{11}$$

which converges for $\rho(L^d D) < 1$.

Consider the residence time $(L + D)^{-1}_{ij}$. The first term on the right-hand side of (11) represents the expected time spent in $i$ according to the regular Markov process, and does not depend on $j$. The next order term $L^d$, however, takes into account both the graph structure and the absorption rates as well as the initial condition $j$. Let $\varepsilon = ||D||/||L||$, and suppose that $\varepsilon \ll 1$. Then, according to (11), $L^d_{ij}$ provides a measure of the deviation in expected residence times between the regular (unperturbed) and transient (perturbed) Markov processes [6]:

$$L^d_{ij} = (\text{Expected time spent in } i \text{ starting from } j) - \frac{1}{d} u_i + \mathcal{O}(\varepsilon). \tag{12}$$

## 3  Properties of $L^d$

In this section we present some basic properties of $L^d$, divided into four parts. First, the definition of $L^d$ and the formula in Proposition 1 involve the spaces $N_{1,0}$ and $R_{1,0}$, which may not be easy to visualize. In Section 3.1 we give a simple formula relating $L^d$ to any $\{1\}$-inverse of $L$. This generalizes a result of Avrachenkov et al. [10], and allows explicit calculation of $L^d$ for several common graph motifs. Second, the interpretation of the entries of $L^d$ in terms of residence time deviations as in (12) is based on the Laurent series expansion (9), and thus is valid only when the absorption rates are small relative to the edge weights of $G$, i.e. $||D|| << ||L||$. In Section 3.2 we show that $L^d$ is related to the fundamental matrix no matter what values the absorptions have relative to the edge weights. Third, in Section 3.3 we provide a topological interpretation for $L^d$ in terms of spanning forests of a related graph, which we call the absorption-scaled graph. These results give a rationale for using $L^d$ to examine the structure of graphs with absorption for arbitrary $d$. Finally, we show that $L^d$ is positive semidefinite for balanced graphs. This will be useful in Section 4 when considering distance metrics based on $L^d$.

First, let us establish that $L^d$ exists and is unique for any strongly connected graph with absorption.

**Theorem 3.** *Let $(G, d)$ be a strongly connected graph with positive absorption vector $d$. Let $L$ denote the graph Laplacian of $G$. Then the absorption inverse $L^d$ exists and is unique.*

*Proof.* Proposition 1 shows existence so let us consider uniqueness. Suppose $X$ and $L^d$ satisfy (1) and (2). We will show that $X$ and $L^d$ agree on $R_{1,0}$ and range $L$, and thus, by Theorem 2b, agree on $\mathbb{R}^n$. Let $v \in R_{1,0}$. By (2), $Xv = 0 = L^d v$, so $X$ and $L^d$ coincide on $R_{1,0}$. Next, let $x \in \text{range } L$. Since $L$ is a bijection from $N_{1,0}$ onto range $L$ (Theorem 2c), there exists unique $y \in N_{1,0}$ such that $Ly = x$. Then we have $Xx = XLy = y$ by (1). On the other hand, $L^d x = L^d Ly = y$, so $L^d$ and $X$ agree on range $L$. Hence $L^d$ and $X$ agree on $\mathbb{R}^n$, and thus $X = L^d$ and $L^d$ is unique. □

### 3.1  Relationship with other $\{1\}$-inverses

In this section we show that $L^d$ can be expressed simply in terms of any $\{1\}$-inverse of $L$. One such matrix that we will consider is the fundamental matrix associated with the regular Markov process generated by $L$. Recall that for a regular Markov chain the fundamental matrix is given by $(I - P + \pi \mathbf{1}^T)^{-1}$ where $P$ is the transition matrix and $\pi$ is the stationary distribution [22]. Similarly, the *fundamental matrix* $Z \in \mathbb{R}^{n \times n}$ for the regular Markov process generated by $L$ is the matrix $Z$ given by

$$Z = (L + \pi w^T)^{-1} = W^{-1}(I - P + \pi \mathbf{1}^T)^{-1}. \tag{13}$$

In Lemma 3 below, we will show that $Z$ is a $\{1\}$-inverse of $L$. To do so, we will need the following result of Boley et al. [2, Lemma 2]:

**Lemma 2** (Boley et al. [2]). *Let $Y$ be a singular matrix, and assume $C = Y + zv^T$ is non-singular. Let $x, y$ be unit vectors such that $Yx = 0$, $Y^T y = 0$. Then, $v^T x \neq 0$, $y^T z \neq 0$, and the inverse of $C$ is*

$$C^{-1} = Y^\dagger - \frac{1}{v^T x} xv^T Y^\dagger - Y^\dagger \frac{1}{y^T z} zy^T + \frac{1 + v^T Y^\dagger z}{v^T x \cdot y^T z} xy^T.$$

**Lemma 3.** *The fundamental matrix, $Z$, for the regular Markov process generated by $L$ satisfies*

$$Z = \left(I - \frac{uw^T}{w^T u}\right) L^\dagger \left(I - \pi \mathbf{1}^T\right) + \frac{u\mathbf{1}^T}{w^T u}. \tag{14}$$

*Furthermore, $Z$ is a $\{1\}$-inverse of $L$.*

*Proof.* The formula (14) follows immediately from Lemma 2 applied to $(L + \pi w^T)^{-1}$. Then, since $Lu = 0$ and $\mathbf{1}^T L = 0$, equation (14) implies

$$LZL = L\left[\left(I - \frac{uw^T}{w^T u}\right) L^\dagger \left(I - \pi \mathbf{1}^T\right) + \frac{u\mathbf{1}^T}{w^T u}\right] L = LL^\dagger L = L.$$

□

The following theorem relates $L^d$ to any $\{1\}$-inverse of the Laplacian matrix. In particular, we can then relate $L^d$ to the familiar matrices $L^\dagger$, $L^\#$, and $Z$. This formula allows for explicit calculation of $L^d$ for some common graph motifs as demonstrated in A.

**Theorem 4.** *Let $Y$ be a $\{1\}$-inverse of the Laplacian matrix, $L$, for a graph with absorption $(G, d)$. Then,*

$$L^d = (I - UD) Y (I - DU). \tag{15}$$

*In particular, we can take $Y = L^\dagger, L^\#,$ or $Z$.*

*Proof.* Using Theorem 2d, Lemma 1 and that $L^d$ is a $\{1, 2\}$-inverse of $L$, we have

$$L^d = L^d L L^d = L^d L Y L L^d = (I - UD) Y (I - DU).$$

We can take $Y = L^\dagger$ or $Y = L^\#$ since they are $\{1, 2\}$-inverses of $L$. In addition, we can take $Y = Z$ by Lemma 3. □

Theorem 4 suggests comparison of $L^d$ with $L^\#$ and $L^\dagger$. Corollary 1 gives necessary and sufficient conditions for equivalence between $L^d$ and the group and Moore-Penrose inverses. For this, recall that a real matrix $H$ is said to be *range Hermitian* if range $H$ = range $H^T$, or equivalently if $\ker H = \ker H^T$ ([17], p. 157). For any matrix $H \in \mathbb{R}^{n \times n}$, $\mathbb{R}^n$ may be written as the following orthogonal direct sums (e.g. Theorem 0.1, [17]):

$$\begin{aligned} \mathbb{R}^n &= \ker H \overset{\perp}{\oplus} \operatorname{range} H^T, \\ \mathbb{R}^n &= \ker H^T \overset{\perp}{\oplus} \operatorname{range} H. \end{aligned} \tag{16}$$

**Corollary 1.** *Let $L^d$ be the absorption inverse for a strongly connected graph with absorption $(G, d)$.*

(a) *$L^d = L^\#$ if and only if the absorption rates are all equal.*

(b) *$L^d = L^\dagger$ if and only if the absorption rates are all equal and $L$ is range Hermitian.*

*Proof.* (a) It is shown in [6] that $L^d = L^\#$ if the absorption rates are equal. If $L^d = L^\#$, then $\ker L = R_{1,0}$ which implies that $D$ is a scalar multiple of $I$. (b) Suppose that $L$ is range Hermitian and the absorption rates are equal. By (a), $L^d = L^\#$, and $L^\# = L^\dagger$ for range Hermitian $L$ (Theorem 4.4, [17]) so $L^d = L^\dagger$.

Now suppose that $L^d = L^\dagger$. It remains to show that this implies equal absorption rates and $L$ range Hermitian. Note that $L^d$ and $L^\dagger$ are both $\{1, 2\}$ inverses, with range $L^d = N_{1,0}$ and $\ker L^d = R_{1,0}$ while range $L^\dagger = \operatorname{range} L^T$ and $\ker L^\dagger = \ker L^T$. By Theorem 2.12 of [17], a $\{1,2\}$ inverse with prescribed range and kernel is unique, and thus $N_{1,0} = \operatorname{range} L^T$ and $R_{1,0} = \ker L^T$.

Let $u$ be a basis for $\ker L$. The columns of $L$ sum to zero, and thus $\ker L^T = \operatorname{Span} \mathbf{1}$. In order to have $R_{1,0} = \ker L^T$, we must then have $d_i = c/u_i$ for some nonzero constant $c$, and thus $D_{ii}^{-1} = c^{-1} u_i$. Consider $N_{1,0}$. For $L^d = L^\dagger$, $N_{1,0} = \operatorname{range} L^T$ and thus $\mathbb{R}^n = \ker L \overset{\perp}{\oplus} D^{-1} \operatorname{range} L$. Let $v \in \operatorname{range} L$. Then $<D^{-1}v, u> = 0$, which gives

$$\sum_{i=1}^n v_i u_i^2 = 0. \tag{17}$$

We now show that this imples $u_i$ the same for all $i$. Let $u^2 = \operatorname{Span}\{(u_1^2, \ldots, u_n^2)^T\}$. Equation (17) implies that $\operatorname{range} L \perp u^2$. Let $x \in u^2$. As $x \perp \operatorname{range} L$, from (16) we have $x \in \ker L^T$. Thus $x$ in $\operatorname{Span} \mathbf{1}$ and $u^2 \subset \operatorname{Span} \mathbf{1}$, giving that the $u_i$ are equal up to sign. But from the Matrix Tree Theorem, all the entries of $u$ have the same sign, and thus all the $u_i$ are equal, so $\ker L = \ker L^T$ and $L$ range Hermitian. Finally, from $d_i = c/u_i$ we have that all the absorption rates are equal. □

We conclude this section by relating $L^d$ to one more matrix that is also related to the graph Laplacian. Let $\hat{L} \in \mathbb{R}^{n-1, n-1}$ denote the matrix obtained by deleting the $n$th row and $n$th column of $L$. The *bottleneck matrix of $L$ based at vertex $n$* is defined as $\hat{L}^{-1}$ [1]. Kirkland et al. have expressed the group inverse of the graph Laplacian in terms of the bottleneck matrix [1, Propositions 2.1, 2.2]. In Proposition 2 below we reformulate their result. This allows us in Corollary 2 to provide a relation between $L^d$ and the bottleneck matrix based at vertex $n$.

**Proposition 2** (Kirkland et al. [1]). *Let $\hat{M} = \hat{L}^{-1}$ be the bottleneck matrix of $L$ based at vertex $n$ and let $M = \begin{bmatrix} \hat{M} & \mathbf{0} \\ \mathbf{0}^T & 0 \end{bmatrix}$. Then,*

$$L^\# = (I - u\mathbf{1}^T) M (I - u\mathbf{1}^T). \tag{18}$$

*Proof.* For a vector $x \in \mathbb{R}^n$, let $\hat{x}$ denote the vector of length $n-1$ given by $(x_1, \ldots, x_{n-1})^T$. It follows from [1, Proposition 2.1] that

$$\begin{aligned} L^\# &= \left(\hat{\mathbf{1}}^T \hat{M} \hat{u}\right) u\mathbf{1}^T + \begin{bmatrix} \hat{M} - \hat{M}\hat{u}\hat{\mathbf{1}}^T - \hat{u}\hat{\mathbf{1}}^T\hat{M} & -\hat{M}\hat{u} \\ -u_n \hat{\mathbf{1}}^T \hat{M} & 0 \end{bmatrix} \\ &= \left(\hat{\mathbf{1}}^T \hat{M} \hat{u}\right) u\mathbf{1}^T + \begin{bmatrix} \hat{M} & \mathbf{0} \\ \mathbf{0}^T & 0 \end{bmatrix} - \begin{bmatrix} \hat{M}\hat{u} \\ 0 \end{bmatrix} \mathbf{1}^T - u \begin{bmatrix} \hat{\mathbf{1}}^T \hat{M} & 0 \end{bmatrix} \\ &= (\mathbf{1}^T M u) u\mathbf{1}^T + M - Mu\mathbf{1}^T - u\mathbf{1}^T M \\ &= (I - u\mathbf{1}^T) M (I - u\mathbf{1}^T). \end{aligned}$$

□

Note that Proposition 2 implies that $M$ is a $\{1\}$-inverse of $L$ since $L = LL^\# L = L(I - u\mathbf{1}^T)M(I - u\mathbf{1}^T)L = LML$. Corollary 2 then follows directly from Theorem 4.

**Corollary 2.** *Let $M$ be defined as in Proposition 2. Then $L^d = (I - UD)M(I - DU)$.*

## 3.2 $L^d$ and the fundamental matrix for arbitrary $d$

The Laurent series given in equation (9) is equivalent to

$$(L + zD)^{-1} = \frac{1}{z} U + (I + zL^dD)^{-1} L^d \tag{19}$$

for $\rho(zL^dD) < 1$. However, $L + zD$ is nonsingular for any $z > 0$ since it has the $Z$ sign pattern and positive column sums, implying that it is a nonsingular $M$ matrix [23]. The propositions that follow here allow us to prove in Theorem 5 that formula (19) actually holds for any $z > 0$, not only in the radius of convergence of the Laurent series as equation (9) does. In the case that $z = 1$, we recover a formula that relates $L^d$ to the fundamental matrix $(L+D)^{-1}$ for arbitrary $d$. This includes, for example, the situation where the absorption rates are large relative to the edge weights of the graph.

**Proposition 3.** *For any $z > 0$, $(L+zD)^{-1}L(I+zL^dD)$ is a projection onto $N_{1,0}$, which sends $\ker L$ to 0. Furthermore,*

$$(L+zD)^{-1}L(I+zL^dD) = I - UD. \tag{20}$$

*Proof.* First, we observe that $L(I+zL^dD)$ sends $\ker L$ to 0 since $L^d$ sends $R_{1,0}$ to 0 by Theorem 2d. Now, let $y \in N_{1,0}$ and $v \in \mathbb{R}^n$ be such that $Dy = Lv$. Then,

$$L(I+zL^dD)y = Ly + zLL^dLv = Ly + zLv = (L+zD)y,$$

which implies that $(L+zD)^{-1}L(I+zL^dD)y = y$. By Theorem 2a, the first statement is proven. Since $I - UD$ is also a projection onto $N_{1,0}$ that sends $\ker L$ to 0 (Theorem 2d), Theorem 2a implies (20). □

**Proposition 4.** *$I + zL^dD$ is nonsingular for any $z > 0$.*

*Proof.* Let $z > 0$. It suffices to show that $x = -zL^dDx$ implies $x = 0$. By Theorem 2a, we write $x = \beta u + y$ for some $\beta \in \mathbb{R}$ and $y \in N_{1,0}$. Let $v \in \mathbb{R}^n$ be such that $Dy = Lv$. Then, since $L^d$ sends $R_{1,0}$ to 0 (Theorem 2d),

$$L^dDx = L^dD(\beta u + y) = L^dDy = L^dLv.$$

Thus, $x = -zL^dDx = -zL^dLv$ implies that $x \in N_{1,0}$ since $L^dL$ is a projection on $N_{1,0}$ by Theorem 2d. Therefore, $\beta = 0$, i.e. $x = y \in N_{1,0}$ and $Dx = Lv$. Further, we have

$$Lx = -zLL^dLv = -zLv = -zDx,$$

i.e. $(L+zD)x = 0$. This implies $x = 0$ since $L + zD$ is nonsingular. □

**Proposition 5.** *For any $z > 0$,*

$$(I + zL^dD)^{-1} = UD + (L+zD)^{-1}L.$$

*Proof.* We first observe that the formula for $L^d$ given in Theorem 4 implies that $DU$ is a right annihilator for $L^d$. This implies $(I+zL^dD)UD = UD$ or, equivalently,

$$UD(I+zL^dD)^{-1} = UD. \tag{21}$$

Right multiplying equation (20) in Proposition 3 by $(I+zL^dD)^{-1}$ gives

$$(L+zD)^{-1}L = (I-UD)(I+zL^dD)^{-1} = (I+zL^dD)^{-1} - UD$$

by (21), so the proof is complete. □

We now proceed with the proof of the main theorem of this section, which shows that $L^d$ is related to the fundamental matrix for any $z > 0$.

**Theorem 5.** *For any $z > 0$,*
$$(L + zD)^{-1} = \frac{1}{z}U + (I + zL^d D)^{-1} L^d.$$

*Proof.* Let $z > 0$. We observe that, by Proposition 5 and since $UDL^d = 0$,
$$\frac{1}{z}U + (I + zL^d D)^{-1} L^d = \frac{1}{z}U + [UD + (L + zD)^{-1} L]L^d = \frac{1}{z}U + (L + zD)^{-1} LL^d.$$

We therefore define $V := \frac{1}{z}U + (L + zD)^{-1} LL^d$ and show that $V = (L + zD)^{-1}$. First, we observe that $(L + zD)V = DU + LL^d = I$ by Lemma 1. It remains to show that $V(L + zD) = I$ for which we consider the terms in $V$ separately. First, $\frac{1}{z}U(L + zD) = UD$ is a projection on ker $L$ by Theorem 2e. Second, we show that $(L + zD)^{-1} LL^d(L + zD)$ is a projection on $N_{1,0}$. Let $y \in N_{1,0}$ and $v \in \mathbb{R}^n$ such that $Dy = Lv$. Then,
$$LL^d(L + zD)y = LL^d Ly + zLL^d Lv = Ly + zLv = (L + zD)y,$$

which implies that $(L + zD)^{-1} LL^d(L + zD)y = y$. This now completes the proof as Theorem 2a then implies that $V(L + zD) = I$. □

## 3.3 The absorption-scaled graph and its spanning forests

In addition to the connections between $L^d$ and the fundamental matrix of the absorbing process (cf. Sections 2.3 and 3.2), there also exist basic connections between the absorption inverse and spanning forests of a modified graph associated with the original graph with absorption $(G, d)$. We refer to the modified graph as the *absorption-scaled graph*.

**Definition 2.** *Let $(G, d)$ be a strongly connected graph with absorption with adjacency matrix $A$. The graph $\tilde{G}$ with adjacency matrix $\tilde{A} = AD^{-1}$ is said to be the absorption-scaled graph associated with $(G, d)$.*

Edge weights of $\tilde{G}$ are thus the edge weights of $G$, rescaled by the absorption rate at the outgoing vertex (i.e. $\tilde{a}_{ij} = a_{ij}/d_j$). The rescaled outdegrees are then given by $\tilde{w} = \tilde{A}^T \mathbf{1} = D^{-1} w$, and we correspondingly define $\tilde{W} = \text{diag}\{\tilde{w}\} = WD^{-1}$. The graph Laplacian of the absorption-scaled graph is $\tilde{L} = LD^{-1}$, and a basis for ker $\tilde{L}$ is given by $\tilde{u} = Du/\bar{d}$.

Proposition 8 below shows that the absorption inverse is closely related to the group inverse of the absorption-scaled graph $\tilde{G}$. We then show in Theorem 6 that $L^d$ has a topological interpretation in terms of the spanning trees and forests of the absorption-scaled graph. These results use the following two propositions, the latter of which gives a new formula for $L^d$.

**Proposition 6.** *For any $z > 0$, $L + zDUD$ is nonsingular.*

*Proof.* Suppose $(L + zDUD)x = 0$. By Theorem 2a, $x = \beta u + y$ for some $\beta \in \mathbb{R}$ and $y \in N_{1,0}$. Then,
$$Ly = Lx = -zDUDx = -z\beta Du \tag{22}$$
where the third equality holds since $UD$ is a projection onto ker $L$ (Theorem 2e). Applying $L^d$ to both sides of (22) gives $y = 0$ by Theorem 2d. It then follows from (22) that $\beta = 0$ and, hence, $x = 0$. □

**Proposition 7.** *For any $z > 0$,*
$$(L + zDUD)^{-1} = \frac{1}{z}U + L^d. \tag{23}$$

*Proof.* Lemma 2 can be used to calculate the inverse since $zDUD = \frac{z}{\bar{d}}Dud^T$ is a rank-one matrix. This gives
$$(L + zDUD)^{-1} = L^\dagger - UDL^\dagger - L^\dagger DU + UDL^\dagger DU + \frac{1}{z}U = (I - UD) L^\dagger (I - DU) + \frac{1}{z}U = L^d + \frac{1}{z}U$$
where the last equality holds by Theorem 4. □

**Proposition 8.** *Let $L^d$ be the absorption inverse for a strongly connected graph with absorption $(G,d)$ and let $\tilde{L}$ be the graph Laplacian of the associated absorption-scaled graph. Then,*

$$\tilde{L}^{\#} = DL^d. \tag{24}$$

*Proof.* We first rewrite Proposition 7 in terms of the absorption-scaled graph:

$$(\tilde{L} + zDU)^{-1} = \frac{1}{z}DU + DL^d. \tag{25}$$

Since $\tilde{u} \in R_{1,0}$ is a basis for $\ker \tilde{L}$, Lemma 1b implies that $DU$ is the eigenprojection at 0 for $\tilde{L}$. Thus, comparing equation (25) and equation (5) gives the result. $\square$

Proposition 8 implies that $\tilde{L}^{\#}$ is similar to $L^d D$ and, therefore, we can express the radius of convergence of the Laurent series (9) in terms of the spectrum of $\tilde{L}^{\#}$. Note also that the relation (24) provides an alternative derivation of the fact that $L^d = L^{\#}$ in the case of equal absorption rates (i.e. Corollary 1a).

The absorption-scaled graph also provides a topological interpretation for the entries of $L^d$ in terms of the forest matrices of $\tilde{G}$. To see this, rewrite the formula given by Theorem 5 in terms of the absorption scaled graph:

$$(z^{-1}\tilde{L} + I)^{-1} = DU + zD(I + zL^dD)^{-1}L^d. \tag{26}$$

By Theorem 1, the left-hand side of (26) can be expressed in terms of forest matrices of $\tilde{G}$. This allows us in Theorem 6 below to express $L^d$ in terms of the matrices of in-forests of $\tilde{G}$ with $n-1$ and $n-2$ arcs, i.e. spanning trees and spanning forests with two trees, respectively. Let $\tilde{Q}_k$ denote the matrix of in-forests with $k$ arcs for the absorption-scaled graph $\tilde{G}$, and $\tilde{J}_k$ denote the normalized matrix of in-forests with $k$ arcs, i.e. $\tilde{J}_k = \tilde{Q}_k/\tilde{\sigma}_k$ where $\tilde{\sigma}_k = \omega(\tilde{\mathcal{F}}_k)$.

**Theorem 6.** *Let $L^d$ be the absorption inverse for the graph with absorption $(G,d)$ and let $\tilde{G}$ be the corresponding absorption-scaled graph. Then,*

$$L^d = \frac{D^{-1}\tilde{\sigma}_{n-2}}{\tilde{\sigma}_{n-1}}\left[\tilde{J}_{n-2} - \tilde{J}_{n-1}\right], \tag{27}$$

*i.e.*

$$L^d_{ij} = \frac{\omega(\tilde{\mathcal{F}}^{j \to i}_{n-2})}{d_i\tilde{\sigma}_{n-1}} - \frac{\tilde{\sigma}_{n-2}u_i}{\tilde{\sigma}_{n-1}\bar{d}}. \tag{28}$$

*Proof.* Let $z > 0$ and $\tau = z^{-1}$. Equation (26) and Theorem 1 (Cheboratev and Agaev [16]) imply

$$\frac{1}{\tilde{\sigma}(\tau)}\left(\tilde{Q}_0 + \tau\tilde{Q}_1 + \ldots + \tau^{n-1}\tilde{Q}_{n-1}\right) = DU + \tau^{-1}D(I + \tau^{-1}L^dD)^{-1}L^d.$$

Left multiplying by $\tau(I + \tau^{-1}L^dD)D^{-1} = \tau D^{-1} + L^d$ gives

$$\frac{1}{\tilde{\sigma}(\tau)}(\tau D^{-1} + L^d)\left(\tilde{Q}_0 + \tau\tilde{Q}_1 + \ldots + \tau^{n-1}\tilde{Q}_{n-1}\right) = \tau U + L^d.$$

Multiplying by $\tilde{\sigma}(\tau) = \sum_{k=0}^{n-1}\tilde{\sigma}_k\tau^k$ gives

$$(\tau D^{-1} + L^d)\left(\tilde{Q}_0 + \tau\tilde{Q}_1 + \ldots + \tau^{n-1}\tilde{Q}_{n-1}\right) = \left(\sum_{k=0}^{n-1}\tilde{\sigma}_k\tau^k\right)(\tau U + L^d)$$

or, equivalently,

$$\sum_{k=0}^{n}\tau^k[L^d\tilde{Q}_k + D^{-1}\tilde{Q}_{k-1}] = \sum_{k=0}^{n}\tau^k\left[\tilde{\sigma}_kL^d + \tilde{\sigma}_{k-1}U\right]$$

where we take $\tilde{Q}_{-1} = \tilde{Q}_n = 0$ and $\tilde{\sigma}_{-1} = \tilde{\sigma}_n = 0$. Equating the coefficients of $\tau^n$ gives

$$\tilde{Q}_{n-1} = \tilde{\sigma}_{n-1}DU. \tag{29}$$

Equating the coefficients of $\tau^{n-1}$ gives

$$L^d \tilde{Q}_{n-1} + D^{-1} \tilde{Q}_{n-2} = \tilde{\sigma}_{n-1} L^d + \tilde{\sigma}_{n-2} U.$$

Note that, by equation (29), $L^d \tilde{Q}_{n-1} = 0$. Then,

$$L^d = \frac{1}{\tilde{\sigma}_{n-1}} \left[ D^{-1} \tilde{Q}_{n-2} - \tilde{\sigma}_{n-2} U \right] = \frac{D^{-1}}{\tilde{\sigma}_{n-1}} \left[ \tilde{Q}_{n-2} - \frac{\tilde{\sigma}_{n-2}}{\tilde{\sigma}_{n-1}} \tilde{Q}_{n-1} \right] = \frac{D^{-1} \tilde{\sigma}_{n-2}}{\tilde{\sigma}_{n-1}} \left[ \tilde{J}_{n-2} - \tilde{J}_{n-1} \right].$$

□

An alternative proof of Theorem 6 follows directly from Proposition 8 and Proposition 15(iii) in [16].

The following simple example demonstrates how the in-forests of $\tilde{G}$ are counted.

**Example** Let $G$ be a symmetric linear path with $n = 3$ vertices and arbitrary absorption $d > 0$ (Figure 1a). The corresponding absorption-scaled graph $\tilde{G}$ is given in Figure 1b. The spanning trees of $\tilde{G}$, i.e. in-forests with $n-1$ arcs, are enumerated with their weights in Figure 1c and the in-forests with $n-2$ arcs in Figure 1d. From the respective weights we compute $\tilde{\sigma}_{n-1} = \frac{d_1 + d_2 + d_3}{d_1 d_2 d_3}$ and $\tilde{\sigma}_{n-2} = \frac{2d_1 d_3 + d_1 d_2 + d_2 d_3}{d_1 d_2 d_3}$. Suppose we wish to calculate $L^d_{12}$. The only in-forest with $n-2$ arcs in which vertex 2 is in a tree converging to vertex 1 is the second forest in Figure 1d. Thus, $\omega(\tilde{\mathcal{F}}^{2 \to 1}_{n-2}) = 1/d_2$ and $L^d_{12}$ can be calculated from equation (28):

$$\begin{aligned}
L^d_{12} &= \frac{1}{\tilde{\sigma}_{n-1}} \left( \frac{\omega(\tilde{\mathcal{F}}^{2 \to 1}_{n-2})}{d_1} - \frac{\tilde{\sigma}_{n-2} u_1}{\bar{d}} \right) \\
&= \frac{d_1 d_2 d_3}{d_1 + d_2 + d_3} \left( \frac{1}{d_1 d_2} - \frac{2d_1 d_3 + d_1 d_2 + d_2 d_3}{d_1 d_2 d_3 n \bar{d}} \right) \\
&= \frac{1}{d_1 + d_2 + d_3} \left( d_3 - \frac{2d_1 d_3 + d_1 d_2 + d_2 d_3}{d_1 + d_2 + d_3} \right) \\
&= \frac{d_3^2 - d_1 d_3 - d_1 d_2}{(d_1 + d_2 + d_3)^2},
\end{aligned}$$

which matches the formulas given in A in Proposition 12 (a path with $n = 3$ can be viewed as a star with vertex 2 as the hub) and Proposition 13.

Finally, we use Theorem 6 to show two more properties of $L^d$ that will be useful in Section 4 for defining a metric based on $L^d$.

**Corollary 3.** $L^d$ has the property of (row) diagonal maximality. That is, for each $i$, $L^d_{ij} < L^d_{ii}$ for all $j \neq i$.

*Proof.* Fix $i$. By Theorem 6, for any $j = 1, ..., n$,

$$L^d_{ij} = \frac{\tilde{\sigma}_{n-2}}{d_i \tilde{\sigma}_{n-1}} \left( \frac{\omega(\tilde{\mathcal{F}}^{j \to *i}_{n-2})}{\tilde{\sigma}_{n-2}} - \frac{d_i u_i}{\bar{d}} \right).$$

The second term is equal for all $j$. The first term is strictly maximized on the diagonal since, for $j \neq i$, $\tilde{\mathcal{F}}^{j \to *i}_{n-2} \subset \tilde{\mathcal{F}}^{i \to *i}_{n-2}$ implies $\omega(\tilde{\mathcal{F}}^{j \to *i}_{n-2}) < \omega(\tilde{\mathcal{F}}^{i \to *i}_{n-2})$. To see that $\tilde{\mathcal{F}}^{j \to *i}_{n-2} \subset \tilde{\mathcal{F}}^{i \to *i}_{n-2}$ is strict, consider a spanning tree converging to $i$. Remove the edge outgoing from $j$. Now we have a forest with two trees with roots $i$ and $j$, and therefore $j$ is not in the tree converging to $i$. Thus, we have found a forest that is in $\tilde{\mathcal{F}}^{i \to *i}_{n-2} \setminus \tilde{\mathcal{F}}^{j \to *i}_{n-2}$. □

**Corollary 4.** $L^d$ has the property of diagonal positivity, i.e. $L^d_{ii} > 0$ for all $i$.

*Proof.* It follows from Theorem 2d that $L^d D u = 0$. Therefore, by Corollary 3, $0 = \sum_j L^d_{ij} d_j u_j < \sum_j L^d_{ii} d_j u_j = L^d_{ii} \bar{d}$. Since $\bar{d} > 0$, the result follows. □

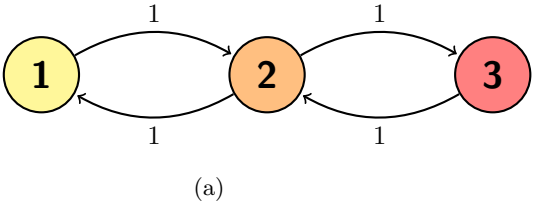
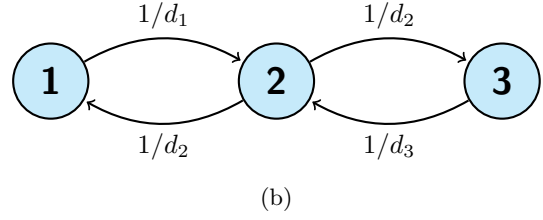
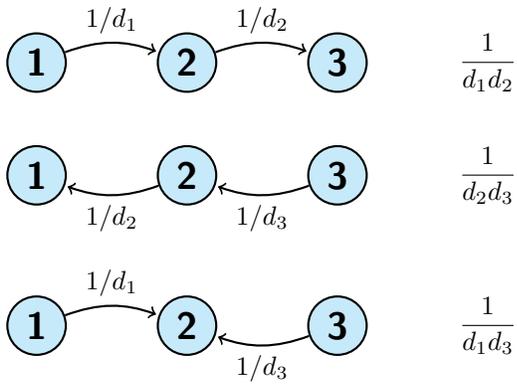
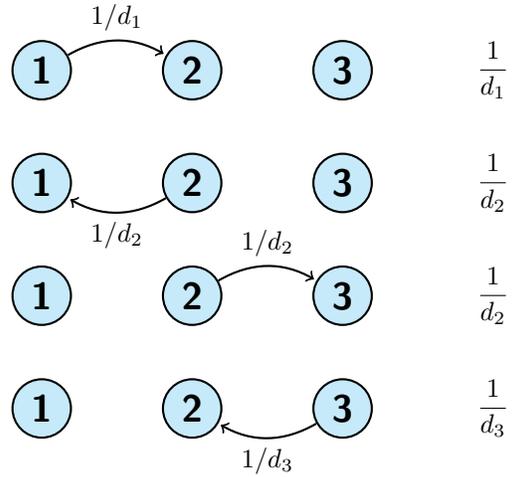

Figure 1: (a) Linear path with $n = 3$ and arbitrary $d > 0$ and (b) its corresponding absorption-scaled graph. Enumeration of the in-forests with (c) $n-1$ and (d) $n-2$ arcs of the absorption-scaled graph, together with their weights.

## 3.4 Balanced graphs and positive semidefiniteness

The final property of $L^d$ included here is that, in the case of a balanced graph, $L^d + (L^d)^T$ is positive semidefinite. This will be useful for relating $L^d$ to a metric for weighted, directed graphs given by Boley et al. [2].

**Lemma 4.** *Let $L$ be the graph Laplacian for a balanced graph. Then, $L + L^T$ is positive semidefinite.*

*Proof.* Consider the Gershgorin discs of $L+L^T$ (e.g. [24, Section 6.1]). Let $\rho_i$ denote the absolute sum of the non-diagonal entries of row $i$ of $L + L^T$. Then,

$$\rho_i = \left| \sum_{j \neq i} (L_{ij} + L_{ji}) \right| = \sum_{j \neq i} a_{ij} + \sum_{j \neq i} a_{ji} = \text{total inflow to } i + \text{total outflow from } i.$$

For balanced graphs, the sum of the inflow equals the sum of the outflow to each vertex, and thus $\rho_i = 2\sum_{j \neq i} a_{ji} = 2L_{ii}$.

The Gershgorin disc theorem then gives that the eigenvalues of $L + L^T$ are contained in the union of discs centered at $2L_{ii}$ with radius $\rho_i$. Therefore, we have that the eigenvalues of $L + L^T$ are contained in $\bigcup_{i=1}^{n}[0, 4\sum_{j \neq i} a_{ji}]$. Thus, the eigenvalues of $L + L^T$ are non-negative, and hence $L + L^T$ is positive semidefinite. □

**Theorem 7.** *Let $L^d$ be the absorption inverse of the graph Laplacian $L$ of a balanced, strongly connected graph. Then $L^d + (L^d)^T$ is positive semidefinite.*

*Proof.* By Boley et al. [2, Lemma 12], this is equivalent to showing that $v^T L^d v \geq 0$ for all $v \in \mathbb{R}^n$. Let $v \in \mathbb{R}^n$. From the decomposition in Theorem 2b, we can write $v = x + y$, where $x \in R_{1,0}$, $y \in \text{range } L$. Note that $L^d x = 0$, giving

$$v^T L^d v = (x+y)^T L^d (x+y) = (x+y)^T L^d y.$$

We will show that $x^T L^d y = 0$ and $y^T L^d y \geq 0$.

Consider first the $x^T L^d y$ term. Since $x \in R_{1,0}$ and the graph is balanced, i.e. $u = (1/n)\mathbf{1}$, we can write $x = \beta Du = \beta d/n$ for some $\beta \in \mathbb{R}$. By Theorem 4, $d^T L^d = 0$ since $d^T(I - UD) = 0$. Thus, $x^T L^d y = 0$.

Now consider the $y^T L^d y$ term. As $L$ is a bijection from $N_{1,0}$ onto range $L$ (Theorem 2c), there exists $b \in N_{1,0}$ such that $Lb = y$. Then we have

$$y^T L^d y = b^T L^T L^d L b = b^T L^T b, \tag{30}$$

as $L^d L$ is a projection on $N_{1,0}$. By Lemma 4, $L + L^T$ is positive semidefinite, which together with Lemma 12 of [2] implies that $b^T L^T b \geq 0$, and hence $v^T L^d v \geq 0$. □

Theorem 7 is relevant for constructing distance metrics on $(G, d)$. For a matrix $Y \in \mathbb{R}^{n \times n}$, let $C_Y : \{1, \ldots, n\} \times \{1, \ldots, n\} \to \mathbb{R}$ be given by

$$C_Y(i, j) = Y_{ii} + Y_{jj} - Y_{ij} - Y_{ji}.$$

As noted by Boley et al. [2, Theorem 16], $C_Y$ is a metric if and only if $Y + Y^T$ is positive semidefinite. In particular, Boley et al. consider the metric $C_{L_0^\dagger}$ where $L_0^\dagger$ is the Moore-Penrose pseudoinverse of what they refer to as the "ordinary Laplacian", $L_0$. In our notation, $L_0 = (I - P)\Pi = LW^{-1}\Pi$; for a balanced graph, $L_0 = L/(\mathbf{1}^T w)$, a scaled version of our graph Laplacian. Boley et al. show that $C_{L_0^\dagger}(i, j)$ is equivalent to the well-known commute time distance between vertices $i$ and $j$ [2, Theorem 15]. Chebotarev et al. also define metrics of this form for multigraphs. Specifically they define "forest distance" as $\frac{1}{2}C_{Q_\tau}$ for $Q_\tau = (I + \tau L)^{-1}$ and "adjusted forest distance" as $\tau C_{Q_\tau}$, which is equivalent to commute time for $\tau \to \infty$ [25].

Since $L^d + (L^d)^T$ is positive semidefinite by Theorem 7, it is natural to consider $C_{L^d}$.

**Proposition 9.** *Let $L^d$ be the absorption inverse of the graph Laplacian $L$ of a balanced, strongly connected graph. Then, $C_{L^d} = \dfrac{1}{\mathbf{1}^T w} C_{L_0^\dagger}$.*

*Proof.* It follows from Theorem 4 that

$$L^d = L^\dagger - L^\dagger DU - UDL^\dagger + UDL^\dagger DU. \tag{31}$$

Note that $[L^\dagger DU]_{ij} = [L^\dagger DU]_{ik}$ for any $j, k$ since $U$ has equal columns. This implies that $C_{L^\dagger DU} = 0$. Likewise, $C_{UDL^\dagger DU} = 0$ and, since the rows of $U$ are also equal for a balanced graph, $C_{UDL^\dagger} = 0$ as well. Equation (31) and the linearity of $C_Y$ in $Y$ then imply that $C_{L^d} = C_{L^\dagger}$. The result then follows from the observation above that $L_0 = L/(\mathbf{1}^T w)$, i.e. $L_0^\dagger = \mathbf{1}^T w L^\dagger$. $\square$

Proposition 9 indicates that for balanced graphs $C_{L^d}$ is equivalent to commute time, a metric which only depends on graph structure. We, on the other hand, seek a notion of distance which takes into account absorption. We present such a metric in Section 4.

## 4 Structure of graphs with absorption

In this section we show how the absorption inverse $L^d$ can be used to describe the structure of balanced graphs with absorption. Specifically, we introduce a distance metric (Section 4.1), a centrality measure (Section 4.2), and a graph partitioning algorithm (Section 4.3) based upon $L^d$. Examples illustrating these measures of graph structure are given in Section 4.4.

### 4.1 Distance metric

Let us begin by defining a directed metric based on $L^d$ for balanced graphs with decay. A *directed metric* (sometimes referred to as a directed semimetric [26] or a quasimetric [27]) on a set $X$ is a nonnegative real-valued function $f : X \times X \to \mathbb{R}_+$ such that $f(x, x) = 0$ for all $x \in X$ and $f(x, y) \leq f(x, z) + f(z, y)$ for all $x, y, z \in X$, i.e. $f$ satisfies the triangle inequality [28]. A symmetry condition is not imposed on $f$.

**Definition 3.** *Let $(G, d)$ be a balanced, strongly connected graph with absorption with absorption inverse $L^d$. Let $K = \max_i L_{ii}^d$. We define $R : V \times V \to \mathbb{R}_+$ to be the function such that*

$$R(j, i) = \begin{cases} K - L_{ij}^d & i \neq j \\ 0 & i = j \end{cases}.$$

In Proposition 10 below we show that $R$ is a directed metric. The proof employs a result of Catral et al. [29]:

**Theorem 8** (Theorem 2.1 [29]). *Let $T \in \mathbb{R}^{n,n}$ be nonnegative, irreducible and stochastic and let $\pi \in \mathbb{R}^n$, $||\pi||_1 = 1$, be the normalized right Perron vector of $T$. Let $B = I - T$. Then, for any $1 \leq i, j, k \leq n$,*

$$\frac{B_{ii}^\#}{\pi_i} - \frac{B_{ij}^\#}{\pi_i} - \frac{B_{ki}^\#}{\pi_k} + \frac{B_{kj}^\#}{\pi_k} \geq 0.$$

**Proposition 10.** *$R$ is a directed metric.*

*Proof.* We observe first that, by Corollaries 3 and 4 and definition of $K$, $R(i, j) > 0$ for $i \neq j$. It remains to show the triangle inequality. Let $1 \leq i, j, k \leq n$. Recall that $\tilde{L} = \tilde{W} - \tilde{A}$ is the Laplacian matrix for the absorption-scaled graph. Let $s = \max_i \tilde{w}_i$. Then, we can rewrite $\tilde{L}$ as

$$\tilde{L} = \tilde{W} - \tilde{A} = sI - \hat{T} = s(I - \tilde{T}) \tag{32}$$

where $\hat{T} = sI - \tilde{W} + \tilde{A}$ and $\tilde{T} = (1/s)\hat{T}$. Note that $sI - \tilde{W} \geq 0$, by definition of $s$, which implies that $\tilde{T}$ is nonnegative. The irreducibility of $\tilde{A}$, which follows from $G$ being strongly connected, implies that $\tilde{T}$ is also irreducible. $\tilde{T}$ is stochastic, i.e. $\mathbf{1}^T \tilde{T} = (1/s)\mathbf{1}^T \hat{T} = (1/s)(s\mathbf{1}^T - \mathbf{1}^T \tilde{L}) = \mathbf{1}^T$. Lastly, recall that $\tilde{u} = Du/\bar{d}$ is a basis for $\ker \tilde{L}$, so $\tilde{T}\tilde{u} = \tilde{u} - (1/s)\tilde{L}\tilde{u} = \tilde{u}$ and $\tilde{u}$ is the right Perron vector of $\tilde{T}$.

Let $\tilde{B} = I - \tilde{T}$ so that $\tilde{L} = s\tilde{B}$ by equation (32). Then, $\tilde{L}^{\#} = (1/s)\tilde{B}^{\#}$ and thus, by Theorem 8,

$$\frac{\tilde{L}_{ii}^{\#}}{\tilde{u}_i} - \frac{\tilde{L}_{ij}^{\#}}{\tilde{u}_i} - \frac{\tilde{L}_{ki}^{\#}}{\tilde{u}_k} + \frac{\tilde{L}_{kj}^{\#}}{\tilde{u}_k} \geq 0. \tag{33}$$

Recall that, by (24), $\tilde{L}_{ij}^{\#} = d_i L_{ij}^d$. Therefore, inequality (33) implies

$$\frac{L_{ii}^d}{u_i} - \frac{L_{ij}^d}{u_i} - \frac{L_{ki}^d}{u_k} + \frac{L_{kj}^d}{u_k} \geq 0. \tag{34}$$

Since $G$ is balanced, $u_i = 1/n$ for all $i$. Thus, (34) implies

$$L_{ij}^d + L_{ki}^d - L_{kj}^d \leq L_{ii}^d \leq K,$$

which implies

$$R(j,k) \leq R(j,i) + R(i,k),$$

i.e. $R$ satisfies the triangle inequality. □

We refer to $R(j,i)$ as the *absorption-scaled forest distance* from $j$ to $i$. As the name suggests, $R$ can be interpreted in terms of the spanning forests of the absorption-scaled graph (by Theorem 6). For a balanced graph, equation (28) in Theorem 6 gives

$$L_{ij}^d = \frac{\omega(\tilde{\mathcal{F}}_{n-2}^{j \to i})}{d_i \tilde{\sigma}_{n-1}} - \frac{\tilde{\sigma}_{n-2}}{\tilde{\sigma}_{n-1} n \bar{d}}. \tag{35}$$

The absorption inverse, and thus $R$, depend upon both the graph structure $G$ as well as the absorption rates $d$. To see the dependence on $G$, note that if vertex $j$ is in many of the trees which converge to $i$ in the in-forests with $n-2$ arcs, we expect $\omega(\tilde{\mathcal{F}}_{n-2}^{j \to i})$ to be large. This corresponds to $L_{ij}^d$ being large and $i$ being close to $j$ (i.e. $R(j,i)$ is small), as expected. To illustrate the effect of absorption on the magnitude of $L_{ij}^d$, consider the simple example of a linear path with 3 vertices, as given in Figure 1a. Consider the distances $R(1,2)$ and $R(3,2)$. As shown in Figure 1d, $\omega(\tilde{\mathcal{F}}_{n-2}^{1 \to 2}) = 1/d_1$ and $\omega(\tilde{\mathcal{F}}_{n-2}^{3 \to 2}) = 1/d_3$. It follows from equation (35) that

$$R(1,2) \geq R(3,2) \quad \Leftrightarrow \quad L_{21}^d \leq L_{23}^d \quad \Leftrightarrow \quad \omega(\tilde{\mathcal{F}}_{n-2}^{1 \to 2}) \leq \omega(\tilde{\mathcal{F}}_{n-2}^{3 \to 2}) \quad \Leftrightarrow \quad d_1 \geq d_3.$$

The absorption-scaled forest distance thus reflects the effect of absorption as expected: if the absorption rate in vertex 1 is higher than that of vertex 3, vertex 1 is considered to be "farther" from the central vertex 2 than is vertex 3.

$R(j,i)$ can also be interpreted in terms of the residence times for the transient random walk on $(G,d)$ (Section 2.3). Consider a transient random walk starting from $j$, in the case where absorption is slow relative to transitions between vertices, i.e. $||D|| \ll ||L||$. For balanced graphs, $u_i = 1/n$, and (12) gives

$$\text{Expected time spent in } i \text{ starting from } j = L_{ij}^d + \frac{1}{n\bar{d}} + \mathcal{O}(z), \tag{36}$$

where $z = ||D||/||L|| \ll 1$. $L_{ij}^d$ thus gives an indication of the time spent in vertex $i$ starting from vertex $j$. Small $R(j,i)$ (i.e. large positive $L_{ij}^d$) corresponds to $i$ and $j$ close, as more time is spent in vertex $i$ starting from $j$ than is predicted by the stationary distribution of the regular Markov process and the average absorption rate. Similarly, large $R(j,i)$ (large negative $L_{ij}^d$) corresponds to $i$ being "far" from $j$, i.e. less time than expected is spent in vertex $i$ starting from $j$.

## 4.2 Centrality measure

Centrality measures quantifying the importance of different vertices are another important class of measures of graph structure [30]. Many different centrality measures exist, including measures based on distance metrics (closeness centrality), geodesics (betweenness centrality), the spectrum of the adjacency matrix (eigenvector centrality), and more.

PageRank is a widely-used centrality measure based on random walks on the graph. As introduced by Brin and Page [12], consider a Markov process on a graph given by two parts: i) a random walk, with transitions corresponding to edges in the graph (e.g. crawling along links between web pages) and ii) "teleportation", where the process jumps according to some probability distribution to different vertices, ignoring the graph structure. The stationary distribution of this Markov process corresponds to the PageRank vector, and gives a measure of the importance (e.g. centrality) of each vertex in the graph. This concept has been extensively used both in practice by search engines such as Google for ranking web pages, and as a general measure of centrality in many other contexts [31–34]. There are deep connections between PageRank and generalized inverses of graph Laplacians. In particular, Chung [35] has explored relationships between PageRank and the group inverse of the graph Laplacian.

The absorption inverse $L^d$ provides a natural extension of PageRank to random walks on graphs with absorption. Specifically, in this section we give an extension of PageRank that is based upon the *quasi-stationary distribution* of the transient random walk generated by a graph with absorption. Similar ideas have been considered by Avrachenkov et al. [10] to devise centrality measures similar to PageRank that are independent of the teleportation probability of the random walk.

**Definition 4.** *Let $(G, d)$ be a balanced, strongly connected graph with absorption. The $L^d$-PageRank vector of $(G, d)$ is $L^d \mathbf{1}$, and the $L^d$-PageRank of vertex $i$ is $\sum_{j=1}^n L_{ij}^d$.*

Note that $L^d$-PageRank provides a relative ranking of centrality of the vertices, even though it is possible for $\sum_{j=1}^n L_{ij}^d$ to be negative for some $i$. This definition of $L^d$-PageRank can be interpreted both in terms of the random walk on $(G, d)$, as well as in terms of the spanning forests as given in Theorem 6. From (36), the $L^d$-PageRank of vertex $i$ is, up to a constant, approximately the time spent in $i$, summed across all initial conditions $j$. Thus, Definition 4 converges to the quasi-stationary distribution in ([10], Definition 1) as $||d|| \to 0$. From Theorem 6, we expect that vertices with low absorption rates $d_i$ and large total weight for in-forests rooted at $i$ to be highly central.

We may also interpret Definition 4 as an approximation to the Perron vector of the generator for the transient random walk on $(G, d)$. To see this, consider the generator in scaled time $\tau = \alpha t$, where $1/\alpha$ is a characteristic time scale for transitions between vertices according to the arc weights of $G$. Specifically, let $\alpha = ||L||$, $\bar{L} = L/\alpha$, $\bar{\varepsilon} = \bar{d}/\alpha$, and $\bar{D} = \text{diag}\{d/\bar{d}\}$. Then the generator becomes

$$L + D = \alpha \bar{L} + \bar{d}\bar{D} = \alpha(\bar{L} + \bar{\varepsilon}\bar{D}),$$

which in scaled time is equal to $\bar{L} + \bar{\varepsilon}\bar{D}$. The generator thus is a perturbation of the Laplacian $\bar{L}$. Consider the case where $\bar{\varepsilon} \ll 1$. Then the zero eigenvalue of $\bar{L}$ perturbs to

$$\lambda = a_0 + a_1 \bar{\varepsilon} + \mathcal{O}(\bar{\varepsilon}^2), \tag{37}$$

with $a_0 = 0$ and $a_1$ computed from the the left eigenvector $\mathbf{1} = (1, \ldots, 1)$ and right eigenvector $u$ of $\bar{L}$ (cf. Horn and Johnson [24], p. 372):

$$\begin{aligned} a_1 &= \tfrac{1}{\langle \mathbf{1}, u \rangle} \sum_{j=1}^n u_j \bar{D}_{jj} \\ &= \sum_{j=1}^n u_j \tfrac{d_j}{\bar{d}} \\ &= 1. \end{aligned} \tag{38}$$

To find the associated eigenvector $p$, look for an eigenvector of the form

$$p \approx cu + \bar{\varepsilon} y + \mathcal{O}(\bar{\varepsilon}^2), \tag{39}$$

where $c > 0$ and $y \in N_{1,0}$. Then

$$(\bar{L} + \bar{\varepsilon}\bar{D})(cu + \bar{\varepsilon} y) = \lambda(cu + \bar{\varepsilon} y) + \mathcal{O}(\bar{\varepsilon}^2). \tag{40}$$

Expanding (40), using (37)−(38) and gathering together $\mathcal{O}(\bar{\varepsilon})$ terms gives

$$\bar{L}y + c\bar{D}u = cu. \tag{41}$$

Multiplying both sides of (41) by $L^d$ on the left and using that $\bar{D}u \in \ker L^d$ gives

$$L^d \bar{L} y = c L^d u. \tag{42}$$

As $L^d \bar{L}$ is a projection on $N_{1,0}$ (Theorem 2d), we have

$$y = c L^d u, \tag{43}$$

and thus
$$p \approx c\left(I + \frac{\bar{d}}{\alpha}L^d\right)u. \tag{44}$$

The stationary distribution of the regular Markov process thus perturbs to the quasi-stationary distribution $p$. For balanced $G$, $u_i = 1/n$ for all $i$, and hence the $p_i$ are distinguished only by the sums of the corresponding rows of $L^d$, giving $L^d$-PageRank centrality as in Definition 4.

## 4.3 Graph partitioning

Another application of $L^d$ is in regards to partitioning algorithms for graphs with absorption $(G, d)$. Consider the classical problem of partitioning an undirected graph without absorption into two groups such that the sum of the edge weights between the groups is small ("min cut") and the two groups are similar in size ("equipartitioning"). As described by von Luxburg [36], the RatioCut introduced by [37] is one way to treat this problem:

$$\min_{H \subset V} \frac{e(H, \bar{H})}{|H||\bar{H}|}, \tag{45}$$

where $H, \bar{H}$ a partition of $G$, $e(H, \bar{H})$ the cut size, and $|H|$ the number of vertices in $H$.

Hagen and Kahng [37] show that this problem can be cast in terms of the Laplacian $L$. Specifically, the minimization problem (45) is equivalent to minimizing the quadratic form $s^T L s$, where (e.g. [36], equation (2))

$$s = \begin{cases} \sqrt{|\bar{H}|/|H|} & v_i \in H \\ -\sqrt{|H|/|\bar{H}|} & v_i \in \bar{H}. \end{cases} \tag{46}$$

A commonly-used approach to solving this NP-hard optimization problem is to instead first solve the "relaxed" problem of minimizing $s^T L s$ over the sphere of radius $\sqrt{n}$ subject to the constraint $s \perp \mathbf{1}$, and then using the resulting solution to produce an indicator vector for group membership [30, 37]. From Rayleigh-Ritz (e.g. [24]), the first minimization step corresponds to computing the eigenvector corresponding to the second smallest eigenvalue of $L$. A crude heuristic for translating this eigenvector into indicators for group membership would be to take group membership as simply the sign of the entries of $s$. This simple method can work well in many cases [38], but refinements for assigning group membership may be preferable for other situations [36].

Note that the group inverse $L^\#$ is a spectral inverse, with the same eigenvectors and corresponding inverse non-zero eigenvalues of $L$ [17]. Thus the relaxed minimization problem has the same solution as the maximization problem $\max_{||s||=\sqrt{n}} s^T L^\# s$.

Now consider a balanced graph with absorption $(G, d)$. The fact that $L^d$ is a generalization of the group inverse suggests an extension of the RatioCut algorithm by replacing $L^\#$ above by $L^d$. Note that

$$s^T L^d s = \frac{1}{2} s^T \left(L^d + (L^d)^T\right) s, \tag{47}$$

and thus by Rayleigh-Ritz maximizing $s^T L^d s$ corresponds to finding the largest eigenvalue of $L^d + (L^d)^T$. Furthermore, $s^T L^d s$ has a natural interpretation in terms of random walks on balanced graphs with absorption. Suppose that $||D|| \ll ||L||$. Then, up to a constant, $L^d_{ij}$ is approximately the residence time spent in $i$ starting from $j$ (cf. equation (12)). Suppose additionally that $|H| = |\bar{H}|$. Then we can take the terms of the membership vector $s$ to be $\pm 1$, giving

$$\begin{aligned} s^T L^d s &= \sum_{i=1}^n \sum_{j=1}^n L^d_{ij} s_i s_j \\ &= \sum_{g_i = g_j} L^d_{ij} - \sum_{g_i \neq g_j} L^d_{ij} \\ &\approx \text{(residence time within communities)} - \text{(residence times between communities)}, \end{aligned} \tag{48}$$

where $g_i, g_j$ denote the groups assigned to vertices $i$ and $j$, respectively.

Setting $s_i = \pm 1$ only corresponds to (46) for equal group sizes. However, as RatioCut (45) favors equal group sizes, good solutions are often obtained by setting $s_i$ to plus or minus one [38]. This then gives the following simple spectral partitioning algorithm for graphs with absorption:

1. Compute $L^d$ from $(G, d)$.

2. Compute the leading eigenvector $s$ of $L^d + (L^d)^T$.

3. Assign vertex $i$ to group 1 if $s_i \geq 0$, and to group 2 if $s_i < 0$.

## 4.4 Examples

This section presents some simple examples to illustrate how the $L^d$-based distance and centrality metrics and partitioning algorithm described in Sections 4.1−4.3 reflect the structure of graphs with absorption. In the figures that follow, absorption rates are indicated by vertex color.

### 4.4.1 Distance

Consider a "bridge" graph in which two three by three grids are connected by a single vertex with equal absorption rates at all vertices (Figure 2a). As discussed in [6], the bridge vertex (labeled 10) acts as a bottleneck to mixing between the left and right sides of the bridge. For equal absorption rates, $L^d = L^\#$ and the absorption-scaled forest distance, $R$, is independent of $d$. We thus expect vertices on opposite sides of the bridge to be far away from each other for a reasonable distance metric. Indeed, $L^d_{ij}$ is positive for $i \neq j$ on the same side of the bridge and is negative for $i$ and $j$ on opposite sides of the bridge, corresponding to small values of $R$ between vertices on the same side of the bridge, and large values for vertices on opposite sides of the bridge (Figure 2b).

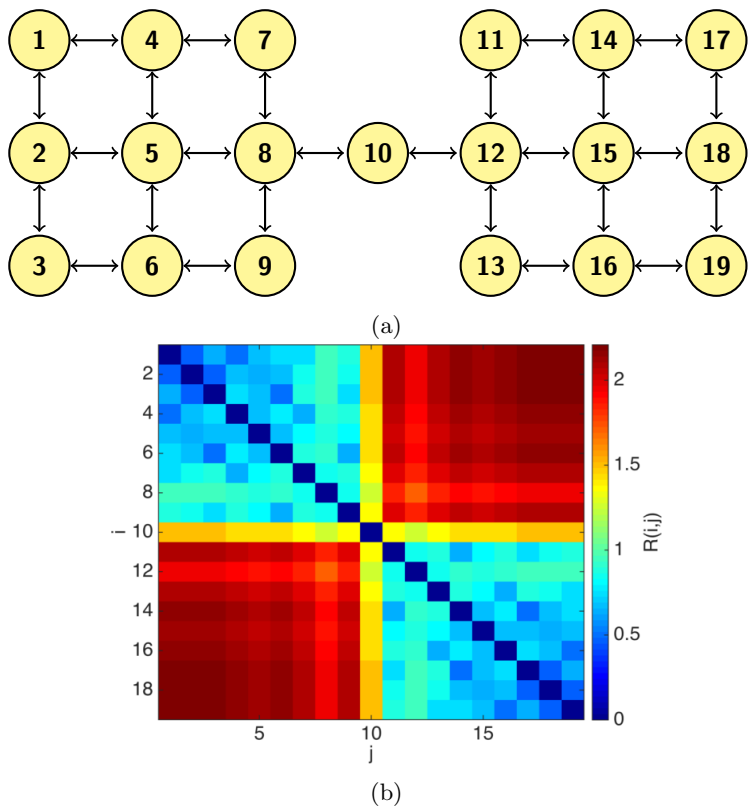

Figure 2: (a) Bridge graph with equal absorption rates ($d = \mathbf{1}$), and (b) the pairwise absorption-scaled forest distances where the color of the $ij$-th cell is determined by $R(j, i)$.

To illustrate the effect of differing absorption rates on $R$, consider next the same bridge graph with a strip of high absorption in vertices 4, 5 and 6 (Figure 3a). Inspection of the pairwise distance matrix (Figure 3b) shows that this strip of high absorption rates separates vertices 1, 2, and 3 from the vertices to the right of the strip. For example, the average absorption-scaled forest distance from vertex 1 to any other vertex in the graph is 2.63. On the other hand, the average distance from vertex 17, which is graphically symmetric to vertex 1, to any other vertex is 1.84. Furthermore, unlike in Figure 2b, we see that vertices 1−9 on the left hand side of the bridge are further from one another (with average distance 2.21) than vertices 11−19 on the right hand side of the bridge (with average distance 0.75).

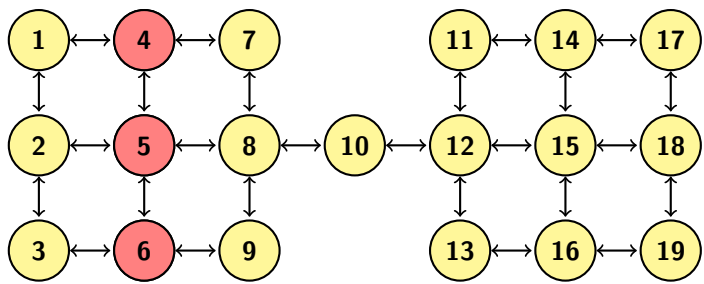

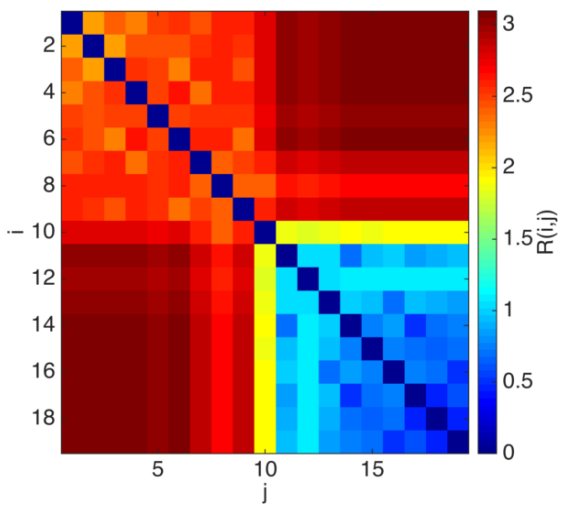

Figure 3: (a) Bridge graph with strip of high absorption, i.e. $d = (1, 1, 1, 10, 10, 10, 1, ..., 1)$, and (b) the pairwise absorption-scaled forest distances where the color of the $ij$-th cell is determined by $R(j, i)$.

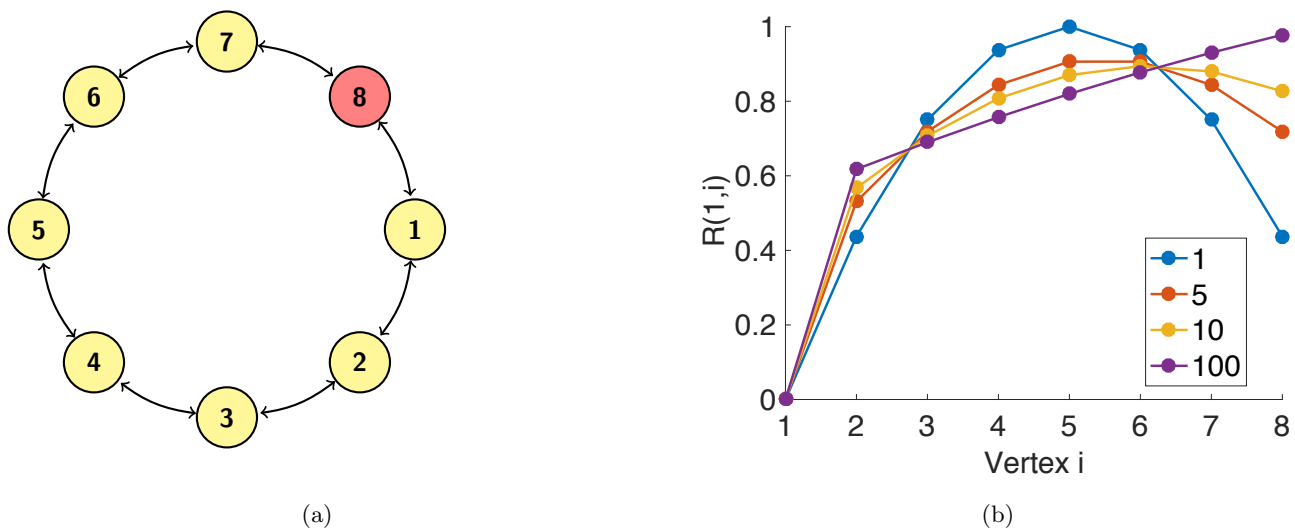

Figure 4: (a) Symmetric cycle with $d_i = 1$ for $i = 1, ..., 7$ and high absorption in vertex 8 and (b) absorption-scaled forest distances, $R(1, i)$, from vertex 1 to vertex $i$ for scenarios of increasing absorption in vertex 8 ($d_8 = 1, 5, 10, 100$ as indicated in legend). Distances are scaled to 1 for each scenario.

The last example we consider is a symmetric cycle with a vertex of high absorption (Figure 4a), as mentioned in the Introduction, to demonstrate how absorption can effectively change the topology of a graph. Consider the absorption-scaled forest distance from vertex 1 to the other vertices as the absorption rate at vertex 8 is varied, with the other absorption rates fixed (Figure 4b). For low absorption in vertex 8, the distance function is non-monotone with vertex number, with vertex 5 being the farthest from 1. However, as $d_8$ increases the distance from 1 becomes monotone, with vertex 8 being the farthest from 1. Thus, the increase in absorption in vertex 8 acts as if to break the cycle and the topology (in terms of absorption-scaled forest distance) effectively becomes a linear arrangement.

### 4.4.2 Centrality

To illustrate how absorption rates can affect vertex centrality, consider a balanced star graph with absorption. Widely used centrality measures such as degree, betweenness, and eigenvector centrality rank the hub (vertex 7 in Figure 5) as most central. However, $L^d$ PageRank can give a different ranking of vertices since it takes absorptions into account in addition to graph structure, which may be important in settings such as disease spread on networks. Figure 5 shows $L^d$-PageRank for two different absorption scenarios for the star. Vertices with a higher relative centrality are indicated with a larger size. Figure 5a shows a star graph where $d_i = 0.1$ for $i = 3, 4, 5, 6, 7$, $d_1 = 1$ and $d_2 = 2$. Computation of $L^d$-PageRank indicates that vertex 2 is the least central and vertex 1 is the second least central, reflecting the high absorption rates at these vertices. In Figure 5b, we take $d_i = 0.1$ for $i = 2, ..., 6$ while $d_1 = d_7 = 0.2$. Note that in this case the hub is more central than vertex 2 despite both vertices having equal absorption rates. This illustrates the dependence of $L^d$-PageRank on both the graph structure as well as the absorption rates. Lastly, note that in the case that all absorptions are equal, $\mathbf{1} \in \ker L^d$ and $L^d$-PageRank ranks all vertices as equally central.

### 4.4.3 Partitioning

Consider spectral partitioning based on $L^d$ (Section 4.3) for a symmetric linear path with absorption, with $n = 8$ vertices and edge weights equal to one. In the first scenario, all absorption rates are taken to be equal, i.e. $d = \mathbf{1}$ (Figure 6a). In the second scenario, the absorption rate in vertex 3 is set to 10 (Figure 6b). Figures 6c and 6d show the components of the dominant eigenvectors of $L^d + (L^d)^T$ in the first and second scenarios, respectively. As expected, when the absorption rates are equal, the graph divides evenly into two groups, $\{1, 2, 3, 4\}$ and $\{5, 6, 7, 8\}$. However, when the absorption rate in vertex 3 is large, spectral partitioning based on $L^d$ results in the groups $\{1, 2, 3\}$ and $\{4, 5, 6, 7, 8\}$. The high absorption rate at 3 increases the distance of vertex 4 from 1 and 2, and thus 4 clusters with $\{5, 6, 7, 8\}$. Numerical inspection indicates a threshold at $d_3 \approx 5.4$ where the partitioning results change from that of the first scenario to that of the second.

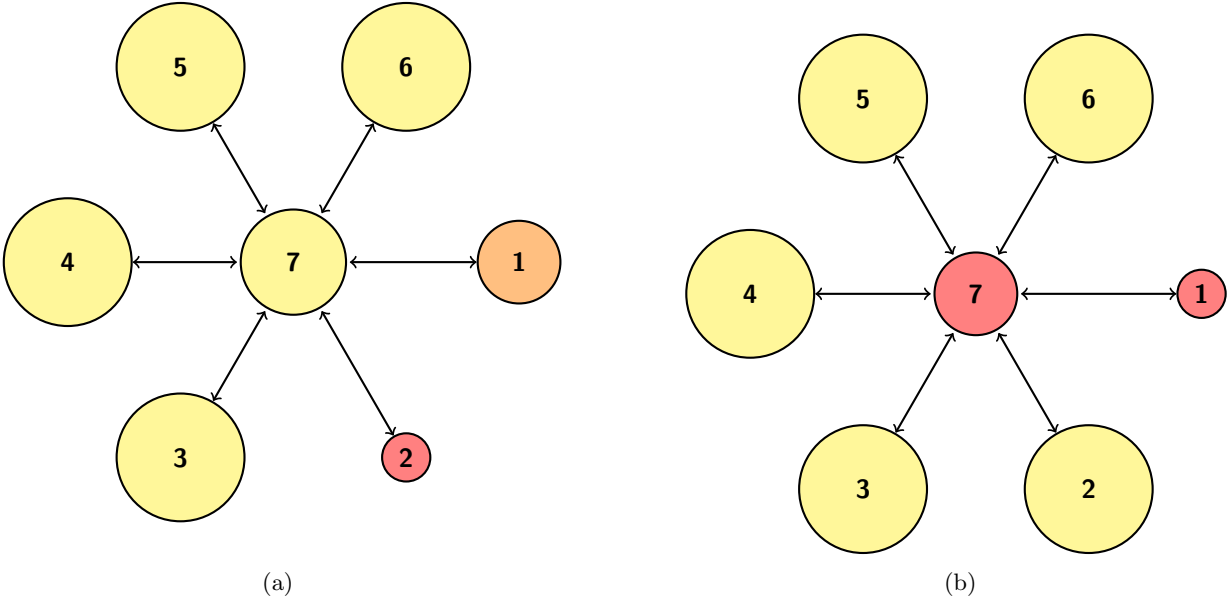

Figure 5: $L^d$-PageRank for a balanced star graph with $n = 7$ and absorption rates given by (a) $d = (1, 2, 0.1, ..., 0.1)$ and (b) $d = (0.2, 0.1, ..., 0.1, 0.2)$. Larger vertices are relatively more central than smaller vertices.

## 5 Discussion

This paper has focused on the absorption inverse of the unnormalized graph Laplacian. Several different Laplacians are widely used, including the "ordinary" Laplacian $L_0 = LW^{-1}\Pi$ [2], and the "symmetric" Laplacian $W^{-1/2}LW^{-1/2}$ [14, 36]. The definition of the absorption inverse (Definition 1) can be correspondingly adjusted for these different Laplacians. For example, let $N_{1,0}^{sym} = W^{1/2}N_{1,0}$, $R_{1,0}^{sym} = W^{-1/2}R_{1,0}$, and $L_{sym}^d = W^{1/2}L^dW^{1/2}$. Then $L_{sym}^d$ is a $\{1, 2\}$ inverse of $L_{sym}$, and indeed is the unique $\{1, 2\}$ inverse of $L_{sym}$ with range $N_{1,0}^{sym}$ and kernel $R_{1,0}^{sym}$ (Theorem 12, chapter 2.6 of [17]). The basic notion of an absorption inverse that incorporates both graph structure as well as absorption rates at the vertices can thus be extended to different Laplacians as appropriate for a particular setting of interest.

The absorption inverse depends upon both graph structure and dynamics (absorption) at the vertices, and descriptors based upon $L^d$ encapsulate both of these features. This combination of graph structure and dynamics is important in many settings [39]. Previous work on finding succinct descriptors of both of these aspects includes work by Ghosh and Lerman [40, 41]. The absorption inverse is promising for understanding the structure of graphs with absorption: $L^d$ arises naturally in the context of the fundamental matrix of the transient random walk on the graph, is connected fundamentally with spanning forests of the absorption-scaled graph, and lends itself easily to several measures of graph structure that take absorption into account, including distance, centrality, and partitioning. Developing efficient algorithms for computing these measures for specific networks and examining their functional significance for processes of interest (for example, contagion processes) is an area for future work. Comparing standard measures of graph structure on the absorption-scaled graph, $\tilde{G}$, with those developed in this work would also be of interest.

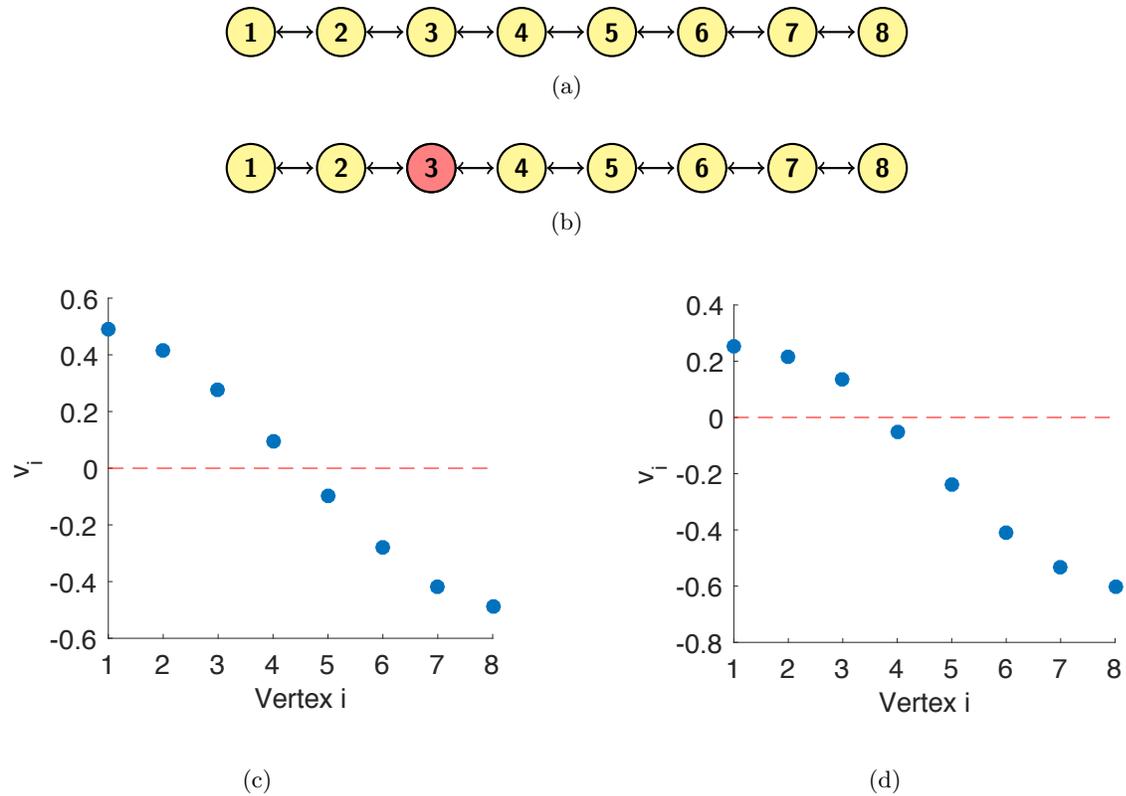

Figure 6: Linear paths with $n = 8$ and absorption rates given by (a) $d = \mathbf{1}$ and (b) $d = (1, 1, 10, 1, ..., 1)$. Entries of the eigenvector $v$ correspond to (c) a partition of graph (a) into $\{1, 2, 3, 4\}$ and $\{5, 6, 7, 8\}$ and (d) a partition of graph (b) into $\{1, 2, 3\}$ and $\{4, 5, 6, 7, 8\}$.

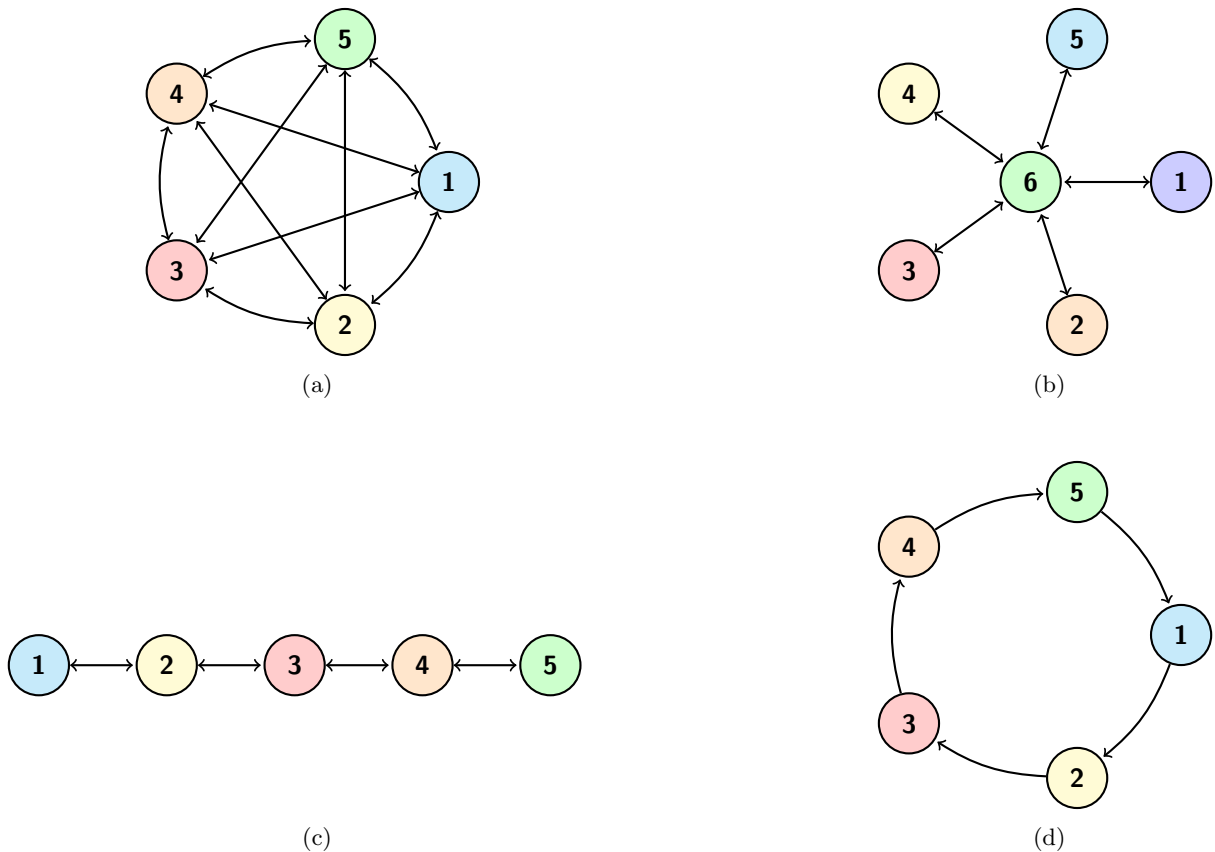

Figure 7: Common graph motifs including (a) complete graph, (b) star graph, (c) linear path, and (d) directed cycle. The multicolored vertices indicate arbitrary absorption rates.

## Acknowledgements

The authors are grateful to Pauline van den Driessche and Mason Porter for helpful discussions. This work was supported by the National Science Foundation (OCE-1115881) and the Mathematical Biosciences Institute (DMS-0931642).

# A  $L^d$ for common graph motifs

By Theorem 4, we are able to explicitly calculate $L^d$ with arbitrary absorption $d > 0$ for a handful of graph motifs: a complete graph where all vertices are pairwise connected, a balanced star, a linear path, and a directed cycle (Figure 7). The following propositions give formulas for $L^d$ in these special cases. Note that the Kronecker delta is denoted by $\delta_{ij}$.

**Proposition 11.** *Let $G$ be a complete graph on $n$ vertices with all edge weights equal to $a$ and let the absorption vector $d > 0$ be arbitrary. Then, the absorption inverse of $L$ with respect to $d$ is given by*

$$L^d = \frac{1}{an}\left(I - \frac{\mathbf{1}d^T}{d^T\mathbf{1}}\right)\left(I - \frac{d\mathbf{1}^T}{d^T\mathbf{1}}\right), \tag{49}$$

*i.e.*

$$L^d_{ij} = \frac{1}{an}\left(\delta_{ij} - \frac{(d_i + d_j)}{d^T\mathbf{1}} + \frac{d^Td}{(d^T\mathbf{1})^2}\right).$$

*Proof.* For such a graph, $L = a(nI - \mathbf{1}\mathbf{1}^T)$ and therefore

$$L + \frac{\mathbf{1}\mathbf{1}^T}{n} = anI + \frac{(1-an)}{n}\mathbf{1}\mathbf{1}^T.$$

By the Sherman-Morrison formula,

$$\left(L + \frac{\mathbf{1}\mathbf{1}^T}{n}\right)^{-1} = \left[anI + \frac{(1-an)}{n}\mathbf{1}\mathbf{1}^T\right]^{-1} = \frac{1}{an}I + \frac{an-1}{an^2}\mathbf{1}\mathbf{1}^T \tag{50}$$

from which it follows from (5) that

$$L^\# = \left(L + \frac{\mathbf{1}\mathbf{1}^T}{n}\right)^{-1} - \frac{\mathbf{1}\mathbf{1}^T}{n} = \frac{1}{an}I - \frac{1}{an^2}\mathbf{1}\mathbf{1}^T. \tag{51}$$

Equation (49) then follows from application of Theorem 4. □

We introduce some notation that will be needed in the proposition below. For a vector $x = (x_1,...,x_n)^T$, let $\hat{x} = (x_1,...,x_{n-1})^T$. Let $\hat{I}$ denote the $n-1$ dimensional identity matrix.

**Proposition 12.** *Let $G$ be a star with $n-1$ leaves, hub vertex labeled $n$, and all edge weights equal to $a$. Let the absorption vector $d > 0$ be arbitrary. Then, the absorption inverse of $L$ with respect to $d$ is given by*

$$L^d = \frac{\hat{d}^T\hat{d}}{a(d^T\mathbf{1})^2}\mathbf{1}\mathbf{1}^T + \frac{1}{ad^T\mathbf{1}}\begin{bmatrix} d^T\mathbf{1}\hat{I} - (\hat{d}\hat{\mathbf{1}}^T + \hat{\mathbf{1}}\hat{d}^T) & -\hat{d} \\ -\hat{d}^T & 0 \end{bmatrix}$$

*i.e.,*

$$a(d^T\mathbf{1})^2 L^d_{ij} = \begin{cases} \hat{d}^T\hat{d} + (d^T\mathbf{1})^2 - 2d^T\mathbf{1}d_i & i = j = \text{leaf} \\ \hat{d}^T\hat{d} - (d^T\mathbf{1})(d_i + d_j) & i,j = \text{leaf}, i \neq j \\ \hat{d}^T\hat{d} - (d^T\mathbf{1})d_j & i = \text{hub}, j = \text{leaf} \\ \hat{d}^T\hat{d} & i = j = \text{hub}. \end{cases}$$

*Proof.* The Laplacian matrix for $G$ is given by

$$L = a\begin{bmatrix} \hat{I} & -\hat{\mathbf{1}} \\ -\hat{\mathbf{1}}^T & n-1 \end{bmatrix}. \tag{52}$$

The bottleneck matrix of $L$ based at vertex $n$, as defined in Section 3.1, is therefore $\hat{M} = \frac{1}{a}\hat{I}$. By Corollary 2,

$$L^d = \frac{1}{a}(I - UD)\begin{bmatrix} \hat{I} & \mathbf{0} \\ \mathbf{0}^T & 0 \end{bmatrix}(I - DU) \tag{53}$$

$$= \frac{1}{a(d^T\mathbf{1})^2}\begin{bmatrix} d^T\mathbf{1}\hat{I} - \hat{\mathbf{1}}\hat{d}^T & -d_n\hat{\mathbf{1}} \\ -\hat{d}^T & d^T\mathbf{1} - d_n \end{bmatrix}\begin{bmatrix} \hat{I} & \mathbf{0} \\ \mathbf{0}^T & 0 \end{bmatrix}\begin{bmatrix} d^T\mathbf{1}\hat{I} - \hat{d}\hat{\mathbf{1}}^T & -\hat{d} \\ -d_n\hat{\mathbf{1}}^T & d^T\mathbf{1} - d_n \end{bmatrix} \tag{54}$$

$$= \frac{1}{a(d^T\mathbf{1})^2}\begin{bmatrix} d^T\mathbf{1}\hat{I} - \hat{\mathbf{1}}\hat{d}^T & \hat{\mathbf{0}} \\ -\hat{d}^T & 0 \end{bmatrix}\begin{bmatrix} d^T\mathbf{1}\hat{I} - \hat{d}\hat{\mathbf{1}}^T & -\hat{d} \\ -d_n\hat{\mathbf{1}}^T & d^T\mathbf{1} - d_n \end{bmatrix} \tag{55}$$

$$= \frac{1}{a(d^T\mathbf{1})^2}\begin{bmatrix} (d^T\mathbf{1})^2\hat{I} - (d^T\mathbf{1})\hat{\mathbf{1}}\hat{d}^T - (d^T\mathbf{1})\hat{d}\hat{\mathbf{1}}^T + \hat{d}^T\hat{d}\hat{\mathbf{1}}\hat{\mathbf{1}}^T & -(d^T\mathbf{1})\hat{d} + \hat{d}^T\hat{d}\hat{\mathbf{1}} \\ -(d^T\mathbf{1})\hat{d}^T + \hat{d}^T\hat{d}\hat{\mathbf{1}}^T & \hat{d}^T\hat{d} \end{bmatrix} \tag{56}$$

from which the result follows. □

The next two propositions use the bottleneck matrix of $L$ based at vertex $n$, $\hat{M}$ as defined in Section 3.1, as well as the matrix $M$ defined in Proposition 2.

**Proposition 13.** *Let $G$ be a linear path with $n$ vertices, all edge weights equal to $a$ and arbitrary absorption vector $d > 0$. Then, the absorption inverse of $L$ with respect to $d$ is given by*

$$L_{ij}^d = n - \max\{i,j\} - \frac{(n-i)}{d^T\mathbf{1}}\sum_{k=1}^{i} d_k - \frac{1}{d^T\mathbf{1}}\sum_{k=i+1}^{n}(n-k)d_k - \frac{(n-j)}{d^T\mathbf{1}}\sum_{k=1}^{j} d_k - \frac{1}{d^T\mathbf{1}}\sum_{k=j+1}^{n}(n-k)d_k$$

$$+ \frac{1}{(d^T\mathbf{1})^2}\sum_{l=1}^{n-1}\left(\sum_{k=1}^{l} d_k\right)^2. \tag{57}$$

*Proof.* The Laplacian matrix for $G$ is given by

$$L = a \begin{bmatrix} 1 & -1 & 0 & \ldots & \ldots & 0 \\ -1 & 2 & -1 & 0 & \ldots & 0 \\ 0 & -1 & 2 & -1 & \ldots & 0 \\ \vdots & \vdots & \vdots & \vdots & \vdots & \vdots \\ 0 & \ldots & \ldots & -1 & 2 & -1 \\ 0 & \ldots & \ldots & \ldots & -1 & 1 \end{bmatrix}.$$

By inspection, we see that

$$\hat{M} = \frac{1}{a}\begin{bmatrix} n-1 & n-2 & n-3 & \ldots & 1 \\ n-2 & n-2 & n-3 & \ldots & 1 \\ n-3 & n-3 & n-3 & \ldots & 1 \\ \vdots & \vdots & \vdots & \vdots & \vdots \\ 1 & 1 & 1 & \ldots & 1 \end{bmatrix},$$

i.e. $M_{ij} = \frac{1}{a}(n - \max\{i,j\})$. Therefore, we can write $M$ as

$$M = \frac{1}{a}\sum_{l=1}^{n-1} \mathbf{1}_l \mathbf{1}_l^T \tag{58}$$

where $\mathbf{1}_l = (1,...,1,0,...,0)^T$, the vector with $l$ ones and $n-l$ zeros. By Corollary 2 and (58),

$$L^d = (I - UD)M(I - DU)$$
$$= M - M\frac{d\mathbf{1}^T}{d^T\mathbf{1}} - \frac{\mathbf{1}d^T}{d^T\mathbf{1}}M + \left(\frac{d^T M d}{(d^T\mathbf{1})^2}\right)\mathbf{1}\mathbf{1}^T$$
$$= \frac{1}{a}\sum_{l=1}^{n-1}\left[\mathbf{1}_l\mathbf{1}_l^T - \frac{\sum_{k=1}^{l} d_k}{d^T\mathbf{1}}(\mathbf{1}_l\mathbf{1}^T + \mathbf{1}\mathbf{1}_l^T) + \frac{\left(\sum_{k=1}^{l} d_k\right)^2}{(d^T\mathbf{1})^2}\mathbf{1}\mathbf{1}^T\right].$$

We observe that $[\mathbf{1}_l\mathbf{1}^T]_{ij} = 0$ if $l < i$ and equals 1 if $l \geq i$. Likewise, $[\mathbf{1}\mathbf{1}_l^T]_{ij} = 0$ if $l < j$ and equals 1 if $l \geq j$. Therefore,

$$\left[\sum_{k=1}^{l} d_k(\mathbf{1}_l\mathbf{1}^T + \mathbf{1}\mathbf{1}_l^T)\right]_{ij} = (n-i)\sum_{k=1}^{i} d_k + \sum_{k=i+1}^{n}(n-k)d_k + (n-j)\sum_{k=1}^{j} d_k + \sum_{k=j+1}^{n}(n-k)d_k$$

from which the formula for the $ij$-th entry, (57), follows. □

**Proposition 14.** *Let $G$ be a directed cycle $1 \to 2 \to ... \to n \to 1$ with edge weights equal to $a$ and arbitrary absorption vector $d > 0$. Then, the absorption inverse of $L$ with respect to $d$ is given by*

$$L_{ij}^d = \frac{1}{a(d^T\mathbf{1})^2}\sum_{l=1}^{n-1}\sum_{k=1}^{l}(d^T\mathbf{1}\delta_{il} - d_l)(d^T\mathbf{1}\delta_{jk} - d_k).$$

*Proof.* The Laplacian matrix for $G$ is

$$L = a \begin{bmatrix} 1 & 0 & \cdots & 0 & -1 \\ -1 & 1 & 0 & \ddots & 0 \\ 0 & -1 & 1 & \ddots & \vdots \\ \vdots & \ddots & \ddots & \ddots & 0 \\ 0 & \cdots & 0 & -1 & 1 \end{bmatrix}. \tag{59}$$

By inspection, we see that

$$\hat{M} = \frac{1}{a} \begin{bmatrix} 1 & 0 & \cdots & 0 & 0 \\ 1 & 1 & 0 & \ddots & 0 \\ 1 & 1 & 1 & \ddots & \vdots \\ \vdots & \ddots & \ddots & \ddots & 0 \\ 1 & \cdots & 1 & 1 & 1 \end{bmatrix}, \tag{60}$$

i.e. $M$ has $\frac{1}{a}$ in the lower triangular part and zeros elsewhere. By Corollary 2 we have

$$L^d = \left(I - \frac{\mathbf{1}d^T}{d^T\mathbf{1}}\right) M \left(I - \frac{d\mathbf{1}^T}{d^T\mathbf{1}}\right) = \frac{1}{(d^T\mathbf{1})^2} QMQ^T \tag{61}$$

where $Q = (d^T\mathbf{1})I - \mathbf{1}d^T$. It follows from (61) that

$$\begin{aligned} L^d_{ij} &= \frac{1}{(d^T\mathbf{1})^2} \sum_{k,l=1}^n Q_{il} Q_{jk} M_{lk} \\ &= \frac{1}{(d^T\mathbf{1})^2} \sum_{k,l=1}^n (d^T\mathbf{1}\delta_{il} - d_l)(d^T\mathbf{1}\delta_{jk} - d_k) M_{lk} \\ &= \frac{1}{a(d^T\mathbf{1})^2} \sum_{l=1}^{n-1} \sum_{k=1}^l (d^T\mathbf{1}\delta_{il} - d_l)(d^T\mathbf{1}\delta_{jk} - d_k). \end{aligned}$$

$\square$